\documentclass[12pt]{amsart}
\usepackage[]{amscd,amsfonts,amsmath,amsxtra,amssymb}
\usepackage{graphics}
\usepackage[all]{xy}
\usepackage{hyperref}
\newtheorem{sub}{}[section]
\newtheorem{subsub}{}[sub]
\usepackage[active]{srcltx}
\def\ov#1{\overline{#1}}

\def\coker{\mathop{\rm coker}\nolimits}
\def\Hom{\mathop{\rm Hom}\nolimits}
\def\HHom{\mathop{\mathcal Hom}\nolimits}

\def\Ext{\mathop{\rm Ext}\nolimits}
\def\EExt{\mathop{\mathcal Ext}\nolimits}
\def\Tor{\mathop{\rm Tor}\nolimits}

\def\deg{\mathop{\rm deg}\nolimits}

\def\Deg{\mathop{\rm Deg}\nolimits}
\def\spec{\mathop{\rm spec}\nolimits}
\def\lra{\longrightarrow}

\def\sigg{\mathop{\hbox{$\displaystyle\sum$}}\limits}

\def\hfl#1#2{\smash{\mathop{\ \hbox to 12mm{\rightarrowfill}}
\limits^{\scriptstyle#1}_{\scriptstyle#2} \ }}
\def\hflb#1#2{\smash{\mathop{\hbox to 12mm{\leftarrowfill}}
\limits^{\scriptstyle#1}_{\scriptstyle#2}}}

\def\m#1{{\hbox{$#1$}}}

\def\ot{\otimes}
\def\og{\leavevmode\raise.3ex\hbox{$\scriptscriptstyle\langle\!\langle$}}
\def\fg{\leavevmode\raise.3ex\hbox{$\scriptscriptstyle\,\rangle\!\rangle$}}
\def\span#1{\langle#1\rangle}

\def\nsp{\lbrace 0\rbrace}

\def\Ssect#1#2{\pagebreak[3]\begin{sub}\label{#2}{\sc\small\small
#1}\rm\medskip}

\def\sepsec{\vskip 1.4cm}
\def\sepsub{\vskip 0.7cm}
\def\sepsubsub{\vskip 0.5cm}
\def\sepprop{\vskip 0.5cm}

\def\xmat#1{\[\xymatrix{#1}\]}
\def\flinc{\ar@{^{(}->}}
\def\fleq{\ar@{=}}
\def\flon{\ar@{->>}}
\def\fmaps{\ar@{|-{>}}}
\def\Nligne{\hfil\break}

\def\BS#1{{\boldsymbol{#1}}}

\newcommand{\C}{{\mathbb C}}

\renewcommand{\P}{{\mathbb P}}

\newcommand{\ka}{{\mathcal A}}
\newcommand{\kb}{{\mathcal B}}
\newcommand{\kc}{{\mathcal C}}
\newcommand{\kd}{{\mathcal D}}
\newcommand{\ke}{{\mathcal E}}
\newcommand{\kf}{{\mathcal F}}

\newcommand{\kh}{{\mathcal H}}
\newcommand{\ki}{{\mathcal I}}
\newcommand{\kj}{{\mathcal J}}
\newcommand{\kk}{{\mathcal K}}

\newcommand{\km}{{\mathcal M}}

\newcommand{\ko}{{\mathcal O}}

\newcommand{\kq}{{\mathcal Q}}

\newcommand{\ks}{{\mathcal S}}

\newcommand{\ku}{{\mathcal U}}

\begin{document}

\def\refname{References}
\def\contentsname{Summary}
\def\proofname{Proof}
\def\abstractname{Resume}

\author{Jean--Marc Dr\'{e}zet}
\address{
Institut de Math\'ematiques de Jussieu,
Case 247,
4 place Jussieu,
F-75252 Paris, France}
\email{drezet@math.jussieu.fr}
\title[{Fragmented deformations}] {Fragmented deformations of primitive
multiple curves}

\begin{abstract}
A {\em primitive multiple curve} is a Cohen-Macaulay irreducible projective 
curve $Y$ that can be locally embedded in a smooth surface, and such that 
$Y_{red}$ is smooth. 

This paper studies the deformations of $Y$ to curves with 
smooth irreducible components, when the number of components is maximal (it is 
then the multiplicity $n$ of $Y$). 

We are particularly interested in deformations to $n$ disjoint smooth 
irreducible components, which are called {\em fragmented deformations}. We 
describe them completely.
We give also a characterization of primitive multiple curves having a 
fragmented deformation.
\end{abstract}

\maketitle
\tableofcontents

\vskip 1cm

\section{Introduction}

A {\em primitive multiple curve} is an algebraic variety $Y$ over $\C$ which is
Cohen-Macaulay, such that the induced reduced variety \m{C=Y_{red}} is a
smooth projective irreducible curve, and that every closed point of $Y$ has a
neighborhood that can be embedded in a smooth surface. These curves have been
defined and studied by C.~B\u anic\u a and O.~Forster in \cite{ba_fo}. The
simplest examples are infinitesimal neighborhoods of projective smooth curves
embedded in a smooth surface (but most primitive multiple curves cannot be
globally embedded in smooth surfaces, cf. \cite{ba_ei}, theorem 7.1).

Let $Y$ be a primitive multiple curve with associated reduced curve $C$, and
suppose that \m{Y\not=C}. Let \m{\ki_C} be the ideal sheaf of $C$ in $Y$. The
{\em multiplicity} of $Y$ is the smallest integer $n$ such that \m{\ki_C^n=0}.
We have then a filtration
\[C=C_1\subset C_2\subset\cdots\subset C_{n}=Y \]
where \m{C_i} is the subscheme corresponding to the ideal sheaf \m{\ki_C^i}
and is a primitive multiple curve of multiplicity $i$. The sheaf \ \m{L=
\ki_C/\ki_C^2} \ is a line bundle on $C$, called the {\em line bundle on $C$
associated to $Y$}.

The deformations of double (i.e. of multiplicity 2) primitive multiple curves
(also called {\em ribbons}) to smooth projective curves have been studied in
\cite{gon}. In this paper we are interested in deformations of primitive
multiple curves \m{Y=C_n} of any multiplicity \m{n\geq 2} to reduced curves 
having
exactly $n$ components which are smooth ($n$ is the maximal number of
components of deformations of $Y$). In this case the number of intersection
points of two components is exactly \m{-\deg(L)}. We give some results in the
general case (no assumption on \m{\deg(L)}) and treat more precisely the case
\m{\deg(L)=0}, i.e. deformations of $Y$ to curves having exactly $n$ disjoint
irreducible components.

\sepsub

\begin{sub}Motivation -- \rm
Let \m{\pi:\kc\to S} be a flat projective morphism of algebraic varieties, $P$
a closed point of $S$ such that \m{\pi^{-1}(P)\simeq Y}, \m{\ko_\kc(1)}
a very ample line bundle on $\kc$ and $P$ a polynomial in one variable with
rational coefficients. Let
\[\tau:\km_{\ko_\kc(1)}(P)\lra S\]
be the corresponding relative moduli space of semi-stable sheaves
(parametrizing the semi-stables sheaves on the fibers of $\pi$ with Hilbert
polynomial $P$ with respect to the restriction of \m{\ko_\kc(1)}, cf.
\cite{si}).

We suppose first that there exists a closed point \m{s\in S} such that 
\m{\kc_s} is a smooth projective irreducible curve. Then in general $\tau$ is 
not flat (some other examples on non flat relative moduli spaces are given in 
\cite{in2}). The reason is that the generic structure of torsion free sheaves 
on $Y$ is more complicated than on smooth curves, and some of these sheaves 
cannot be deformed to sheaves on the smooth fibers of $\pi$.

A coherent sheaf on a smooth algebraic variety 
$X$ is locally free on some nonempty open subset of $X$. This is not true on 
$Y$. But a coherent sheaf $E$ on $Y$ is {\em quasi locally free} on some 
nonempty open subset of $Y$, i.e. on this open subset, $E$ is locally 
isomorphic to a sheaf of the form \ \m{\bigoplus_{1\leq i\leq n}m_i\ko_{C_i}}, 
the sequence of non negative integers \m{(m_1,\ldots,m_n)} being uniquely 
determined (cf. \cite{dr2}, \cite{dr4}). It is not hard to see that if $E$ can 
be extended to a coherent sheaf on $\kc$, flat on $S$, then \ 
\m{R(E)=\sigg_{i=1}i.m_i} \ must be a multiple of $n$.
For example,  it is impossible to deform the stable sheaf
\m{\ko_{C_i}} on $Y$ in sheaves on the smooth fibers, if \m{1\leq i<n}. 

Now suppose that all the fibers \m{\pi^{-1}(s)}, \m{s\not=P}, are reduced with 
exactly $n$ smooth components. I conjecture that (with suitable hypotheses) a 
torsion free coherent sheaf on $Y$ can be extended to a coherent sheaf on 
$\kc$, flat on $S$, using the fact that we allow coherent sheaves of the 
reducible fibers \m{\kc_s} that have not the same rank on all the components. 
This would be a step in the study of the flatness of 
$\tau$. For example (for suitable $\pi$), there exists 
a coherent sheaf $\ke$ on $\kc$, flat on $S$, such that \m{\ke_P=\ko_C}, and 
that for \m{s\not=P}, \m{\ke_s} is the structural sheaf of an irreducible 
component of \m{\kc_s}.

Moduli spaces of sheaves on reducible curves have been studied in \cite{tei1}, 
\cite{tei2}, \cite{tei3}.
\end{sub}

\sepsub

\begin{sub}Maximal reducible deformations -- \rm
Let \m{(S,P)} be the germ of a smooth curve. Let $Y$ be a primitive multiple 
curve
of multiplicity \m{n\geq 2} and \m{k>0} an integer. Let \ \m{\pi:\kc\to S} be
a flat morphism, where $\kc$ is a reduced algebraic variety, such that
\begin{enumerate}
\item[--] For every closed point $s\in S$ such that $s\not=P$, the fiber
$\kc_s$ has $k$ irreducible components, which are smooth and transverse, and
any three of these components have no common point.
\item[--] The fiber $\kc_P$ is isomorphic to $Y$.
\end{enumerate}
\end{sub}
We show that by making a change of variable, i.e. by considering a suitable
germ \m{(S',P')} and a non constant morphism \m{\tau:S'\to S}, and replacing
$\pi$ with \m{\tau^*\kc\to S'}, we can suppose that $\kc$ has exactly $k$
irreducible components, inducing on every fiber \m{\kc_s}, \m{s\not=P} the $k$
irreducible components of \m{\kc_s}. In this case $\pi$ is called a {\em
reducible deformation of $Y$ of length $k$}.

We show that \m{k\leq n}. We say that $\pi$ (or $\kc$) is a {\em maximal
reducible deformation of $Y$} if \m{k=n}.

Suppose that $\pi$ is a maximal reducible deformation of $Y$. We show that if
\m{\kc'} is the union of \m{i>0} irreducible components of $\kc$, and
\m{\pi':\kc'\to S} is the restriction of $\pi$, then \m{\pi'^{-1}(P)\simeq
C_i}, and \m{\pi'} is a maximal  reducible deformation of \m{C_i}. Let
\m{s\in S\backslash\{P\}}. We prove that the irreducible components of $\kc_s$
have the same genus as $C$. Moreover, if \m{D_1,D_2} are distinct irreducible
components of \m{\kc_s}, then \m{D_1\cap D_2} consists of \m{-\deg(L)} points.

\sepsub

\begin{sub}Fragmented deformations (definition) -- \rm Let $Y$ be a primitive
multiple curve of multiplicity \m{n\geq 2} and \m{\pi:\kc\to S} a maximal
reducible deformation of $Y$. We call it a {\em fragmented deformation of $Y$}
if \m{\deg(L)=0}, i.e. if for every \m{s\in S\backslash\{P\}}, \m{\kc_s} is the
disjoint union of $n$ smooth curves. In this case $\kc$ has $n$ irreducible
components \m{\kc_1,\ldots,\kc_n} which are smooth surfaces.

The variety $\kc$ appears as a particular case of a {\em gluing} of
\m{\kc_1,\ldots,\kc_n} along $C$ (cf. \ref{ecl2c}). We prove (proposition
\ref{ecl2b}) that such a gluing $\kd$ is a fragmented deformation of a
primitive multiple curve if and only if every closed point in $C$ has a
neighborhood in $\kd$ that can be embedded in a smooth variety of dimension
3. The simplest gluing is the trivial or {\em initial gluing} $\BS{\ka}$. An
open subset $U$ of $\BS{\ka}$ (and $\kc$) is given by open subsets
\m{U_1,\ldots,U_n} of \m{\kc_1,\ldots,\kc_n} respectively, having the same
intersection with $C$, and
\[\ko_\BS{\ka}(U) \ = \
\{(\alpha_1,\ldots,\alpha_n)\in\ko_{\kc_1}(U\cap\kc_1)\times
\cdots\times\ko_{\kc_n}(U\cap\kc_n) ;\quad \alpha_{1\mid 
C}=\cdots=\alpha_{n\mid C}\} ,\]
and \m{\ko_\kc(U)} appears as a subalgebra of \m{\ko_\BS{\ka}(U)}, hence we
have a canonical morphism \m{\BS{\ka}\to\kc}.

We can view elements of \m{\ko_\kc(U)} as $n$-tuples
\m{(\alpha_1,\ldots,\alpha_n)}, with \m{\alpha_i\in\ko_{\kc_i}(U\cap\kc_i)}.
In particular we can write \ \m{\pi=(\pi_1,\ldots,\pi_n)}.

\end{sub}

\sepsub

\begin{sub} A simple analogy -- \rm Consider $n$ copies of $\C$ glued at 0.
Two extreme examples appear: the trivial gluing \m{\BS{\ka}_0} (the set of
coordinate lines in \m{\C^n}), and a set \m{\kc_0} of $n$ lines in $\C^2$. We
can easily construct a bijective morphism \m{\BS{\Psi}:\BS{\ka}_0\to\kc_0}
sending each coordinate line to a line in the plane

\bigskip

\includegraphics{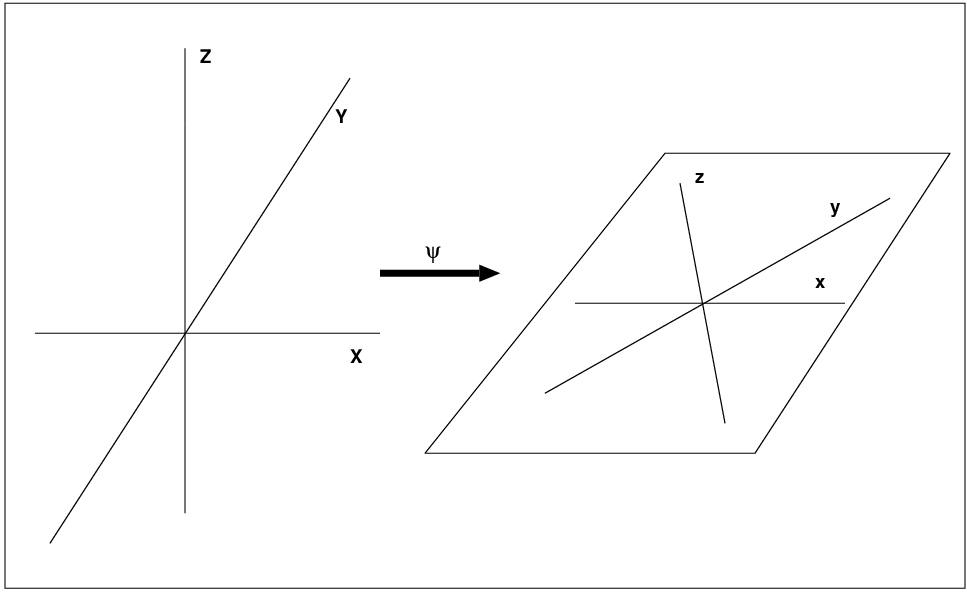}

\bigskip

But the two schemes are of course not isomorphic: the maximal ideal of 0 in
\m{\BS{\ka}_0} needs $n$ generators, but 2 are enough for the maximal ideal of
0 in \m{\kc_0}.

Let \m{\pi_{\kc_0}:\kc_0\to\C} be a morphism sending each component linearly 
onto $\C$, and\Nligne 
\m{\pi_{\BS{\ka}_0}=\pi_{\kc_0}\circ\BS{\Psi}:\BS{\ka}_0\to\C}. 
The difference of \m{\BS{\ka_0}} and \m{\kc_0} can be also seen by using the 
fibers of 0: we have
\[\pi_{\kc_0}^{-1}(0)\simeq\spec(\C[t]/(t^n)) \quad \text{and} \quad
\pi_{\BS{\ka}_0}^{-1}(0)\simeq\spec(\C[t_1,\ldots,t_{n-1}]/(t_1,\ldots,
t_{n-1})^2) \ .\]
Let $\kd$ be a general gluing of $n$ copies of $\C$ at 0, such that there
exists a morphism \m{\pi:\kd\to\C} inducing the identity on each copy of $\C$.
It is easy to see that we have \ \m{\pi^{-1}(0)\simeq\spec(\C[t]/(t^n))} \ if
and only if some neighborhood of 0 in $\kd$ can be embedded in a smooth
surface.
\end{sub}

\sepsub

\begin{sub}Fragmented deformations (main properties) -- \rm Let \m{\pi:\kc\to
S} be a fragmented deformation of \m{Y=C_n}. Let \m{I\subset\{1,\ldots,n\}} be
a proper subset, \m{I^c} its complement, and \m{\kc_I\subset\kc} the subscheme
union of the \m{\kc_i, i\in I}. We prove (theorem \ref{ecl19}) that the ideal
sheaf \m{\ki_{\kc_I}} of \m{\kc_I} is isomorphic to \m{\ko_{\kc_{I^c}}.}

In particular, the ideal sheaf \m{\ki_{\kc_i}} of \m{\kc_i} is generated by a
single regular function on $\kc$. We show that we can find such a generator
such that for \m{1\leq j\leq n}, \m{j\not=i}, its $j$-th coordinate  can be
written as \m{\alpha_j\pi_j^{p_{ij}}}, with \m{p_{ij}>0} and \m{\alpha_j\in 
H^0(\ko_S)} such that \m{\alpha_j(P)\not=0}. If \m{1\leq j\leq n} and 
\m{j\not=i}, we can then obtain a generator that can be written as
\[{\bf u}_{ij} \ = \ (u_1,\ldots,u_m) ,\]
with
\[u_i=0, \quad u_m \ = \ \alpha_{ij}^{(m)}\pi_m^{p_{im}} \ \ \text{for} \
m\not=i, \quad \alpha_{ij}^{(j)}=1 . \]
The constants \ \m{{\bf a}_{ij}^{(m)}=\alpha_{ij\mid C}^{(m)}\in\C} \ have
interesting properties (propositions \ref{ecl17}, \ref{ecl11}). Let
\m{p_{ii}=0} for \m{1\leq i\leq n}. The symmetric matrix \m{(p_{ij})_{1\leq
i,j\leq n}} is called the {\em spectrum} of $\pi$ (or $\kc$).

It follows also from the fact that \m{\ki_{\kc_i}=({\bf u}_{ij})} that $Y$ is
a {\em simple} primitive multiple curve, i.e. the ideal sheaf of $C$ in
\m{Y=C_n} is isomorphic to \m{\ko_{C_{n-1}}}. Conversely, we show in theorem
\ref{theo_sim} that if $Y$ is a simple primitive multiple curve, then there
exists a fragmented deformation of $Y$.

We give in \ref{const_def} and \ref{const_def2} a way to construct fragmented
deformations by induction on $n$. This is used later to prove statements on
fragmented deformations by induction on $n$.
\end{sub}

\sepsub

\begin{sub} $n$-stars and structure of fragmented deformations -- \rm An {\em
$n$-star} of \m{(S,P)} is a gluing $\BS{\ks}$ of $n$ copies \m{S_1,\ldots,S_n} 
of $S$ at $P$, together with a morphism \m{\tau:\BS{\ks}\to S} which is an 
identity on each \m{S_i}. All the $n$-stars have the same underlying 
Zariski topological space \m{S(n)}.

An $n$-star is called {\em oblate} if some neighborhood of $P$ can be
embedded in a smooth surface. This is the case if and only \
\m{\tau^{-1}(0)\simeq\spec(\C[t])/(t^n)}.

Oblate $n$-stars are analogous to fragmented deformations but simpler. We
provide a way to build oblate $n$-stars by induction on $n$.

Let \m{\pi:\kc\to S} be a fragmented deformation of \m{Y=C_n}. We associate to
it an oblate $n$-star $\BS{\ks}$ of $S$. Let \m{\kc^{top}} be the Zariski 
topological space of $\kc$. We have an obvious continuous map \ 
\m{\widetilde{\pi}:\kc^{top}\to S(n)}. For every open subset $U$ of
\m{S(n)}, \m{\ko_\BS{\ks}(U)} is the set of
\m{(\alpha_1,\ldots,\alpha_n)\in\ko_\kc(\widetilde{\pi}^{-1}(U))} such that
\m{\alpha_i\in\ko_{S_i}(U\cap S_i)} for \m{1\leq i\leq n}. We obtain also
a canonical morphism \ \m{\BS{\Pi}:\kc\lra\BS{\ks}} .
We prove (theorem \ref{st_fr_2}) that $\BS{\Pi}$ is flat. Hence it is a flat
family of smooth curves, with \ \m{\BS{\Pi}^{-1}(P)=C}. The converse is also
true, i.e. starting from an oblate $n$-star of $S$ and a flat family of smooth
curves parametrized by it, we obtain a fragmented deformation of a multiple
primitive curve of multiplicity $n$.
\end{sub}

\sepsub

\begin{sub} Fragmented deformations of double curves -- \rm Let \m{Y=C_2} be a
primitive double curve, $C$ its associated smooth curve, \m{\pi:\kc\to S} a
fragmented deformation of $Y$, of spectrum \m{\begin{pmatrix} 0 & p\\
p & 0\end{pmatrix}}, and \m{\kc_1}, \m{\kc_2} the irreducible components of
$\kc$. For \m{i=1,2}, \m{q>0}, let \m{C^q_i} be the infinitesimal
neighborhood of order $q$ of $C$ in \m{\kc_i} (defined by the ideal
sheaf \m{(\pi_i^q)}). It is a primitive multiple curve of multiplicity $q$.

It follows from \ref{ecl12} that \m{C_1^p} and \m{C_2^p} are isomorphic, and
\m{C_1^{p+1}}, \m{C_2^{p+1}} are two extensions of \m{C_1^p} in primitive
multiple curves of multiplicity \m{p+1}. According to \cite{dr1} these
extensions are parametrized by an affine space with associated vector space
\m{H^1(C,T_C)} (where \m{T_C} is the tangent bundle of $C$). Let \m{w\in
H^1(C,T_C)} be the vector from \m{C_1^{p+1}} to \m{C_2^{p+1}}.

Similarly, the primitive double curves with associated smooth curve $C$ such
that \m{\ki_C\simeq\ko_C} are parametrized by \m{\P(H^1(C,T_C))\cup\nsp} (cf.
\cite{ba_ei}, \cite{dr1}).

We prove in theorem \ref{theo_cf} that the point of \m{\P(H^1(C,T_C))\cup\nsp}
corresponding to \m{C_2} is \m{\C w}.
\end{sub}

\sepsub

\begin{sub}\label{notat} Notation: \rm Let $X$ be an algebraic variety and
\m{Y\subset X} a closed subvariety. We will denote by \m{\ki_{Y,X}} (or
\m{\ki_Y} if there is no risk of confusion) the ideal sheaf of $Y$ in $X$.
\end{sub}

\sepsec

\section{Preliminaries}\label{prelim}

\Ssect{Local embeddings in smooth varieties}{plongt}

\begin{subsub}\label{plongt1}{\bf Proposition: } Let $X$ be an algebraic
variety, $x$ a closed point of $X$ and $n$ a positive integer. Then the three
following properties are equivalent:
\begin{enumerate}
\item[(i)] There exist a neighborhood $U$ of $x$ and an embedding $U\subset
Z$ in a smooth variety of dimension $n$.
\item[(ii)] The $\ko_{X,x}$-module $m_{X,x}$ (maximal ideal of $x$) can be
generated by $n$ elements.
\item[(iii)] We have \ $\dim_\C(m_{X,x}/m_{X,x}^2)\leq n$.
\end{enumerate}
\end{subsub}

\begin{proof} It is obvious that (i) implies (ii), and (ii),(iii) are
equivalent according to Nakayama's lemma. It remains to prove that
(iii) implies (i).

Suppose that (iii) is true. There exist an integer $N$ and an embedding \
\m{X\subset\P_N}. Let \m{\ki_X} be the ideal sheaf of $X$ in \m{\P_N}. Let $p$
be the biggest integer such that there exists \m{f_1,\cdots,f_p\in
\ki_{X,x}} whose images in the $\C$-vector space \
\m{m_{\P_N,x}/m_{\P_N,x}^2} \ are linearly independent. Then we have
\[\ki_{X,x} \ \subset \ (f_1,\cdots,f_p)+m_{\P_N,x}^2 .\]
In fact, let \ \m{f\in\ki_{X,x}}. Since $p$ is maximal, the image of $f$
in \ \m{m_{\P_N,x}/m_{\P_N,x}^2} \ is a linear combination of those of
\m{f_1,\cdots,f_n}. Hence we can write
\[f=\sigg_{i=1}^p\lambda_if_i+g , \quad \text{with} \quad \lambda_i\in\C, \
g\in m_{\P_N,x}^2 \ ,\]
and our assertion is proved. It follows that we have a surjective morphism
\[\alpha:\ko_{X,x}/m_{X,x}^2\lra\ko_{\P_N,x}/\big((f_1,\cdots,f_p)+
m_{\P_N,x}^2\big) \ .\]
We have
\[\dim_\C(\ko_{X,x}/m_{X,x}^2) \ \leq \ n+1 , \quad
\dim_\C\big(\ko_{\P_N,x}/\big((f_1,\cdots,f_p)+m_{\P_N,x}^2\big)\big) \ =
N-p+1 \ .\]
Hence \ \m{N-p+1\leq n+1}, i.e. \ \m{p\geq N-n}. We can take for $Z$ a
neighborhood of $x$ in the subvariety of \m{\P_N} defined by
\m{f_1,\cdots,f_{N-n}}, which is smooth at $x$.
\end{proof}

\end{sub}

\sepsub

\Ssect{Flat families of coherent sheaves}{cohpla}

Let \m{(S,P)} be the germ of a smooth curve and \m{t\in\ko_{S,P}} a generator 
of the maximal ideal. Let \ \m{\pi:X\to S} \ be a flat morphism. If $\ke$ is a
coherent sheaf on $X$, $\ke$ is flat on $S$ at \m{x\in\pi^{-1}(P)} if and only
if the multiplication by \m{t:\ke_x\to \ke_x} is injective. In particular the
multiplication by \m{t:\ko_{X,x}\to\ko_{X,x}} is injective.

\sepprop

\begin{subsub}\label{defcmp2}{\bf Lemma: } Let $\ke$ be a coherent sheaf
on $X$ flat on $S$. Then, for every open subset $U$ of $X$, the restriction \
\m{\ke(U)\to\ke(U\backslash\pi^{-1}(P))} \ is injective.
\end{subsub}

\begin{proof} Let \m{s\in\ke(U)} whose restriction to
\m{U\backslash\pi^{-1}(P)} vanishes. We must show that \m{s=0}. By covering
$U$ with smaller open subsets we can suppose that $U$ is affine:
\m{U=\spec(A)}. Hence \ \m{U\backslash\pi^{-1}(P)=\spec(A_t)}. Let
\m{M=\ke(U)}, it is an $A$-module. We have \ \m{\ke_{\mid U}=\widetilde{M}} and
\m{\ke(U\backslash\pi^{-1}(P))=M_t}. Hence if the restriction of $s$ to
\m{U\backslash\pi^{-1}(P)} vanishes, there exists an integer \m{n>0} such that
\m{t^ns=0}. Since the multiplication by $t$ is injective (because $\ke$ is flat
on $S$), we have \m{s=0}. \end{proof}

\sepprop

Let $\ke$ be a coherent sheaf on $X$ flat on $S$. Let \ \m{\kf\subset\ke_{\mid
X\backslash\pi^{-1}(P)}} \ be a subsheaf. For every open subset $U$ of $X$ we
denote by \m{\overline{\kf}(U)} the subset of \m{\kf(U\backslash\pi^{-1}(P))}
of elements that can be extended to sections of $\ke$ on $U$. If \m{V\subset
U} is an open subset, the restriction \
\m{\kf(U\backslash\pi^{-1}(P))\to\kf(V\backslash\pi^{-1}(P))} \ induces a
morphism \ \m{\overline{\kf}(U)\to\overline{\kf}(V)} .

\sepprop

\begin{subsub}\label{defcmp1}{\bf Proposition: } $\overline{\kf}$ is a
subsheaf of $\ke$, and \m{\ke/\overline{\kf}} is flat on  $S$.
\end{subsub}

\begin{proof}
To prove the first assertion, we must show that if $U$ is an open subset of
$X$ and \m{(U_i)_{i\in I}} is an open cover of $U$, then
\begin{enumerate}
\item[(i)] If $s\in\overline{\kf}(U)$ is such that for every $i$ we have
$s_{\mid U_i}=0$, then $s=0$.
\item[(ii)] For every $i\in I$ let $s_i\in\overline{\kf}(U_i)$. Then if for
all $i$, $j$ we have \ $s_{i\mid U_{ij}}=s_{j\mid U_{ij}}$ , then there exists
$s\in\overline{\kf}(U)$ such that for every $i\in I$ we have \ $s_{\mid
U_i}=s_i$.
\end{enumerate}
This follows easily from lemma \ref{defcmp2}.

Now we prove that \m{\ke/\overline{\kf}} is flat on $S$. Let
\m{x\in\pi^{-1}(P)} and \m{u\in(\ke/\overline{\kf})_x} such that \m{tu=0}. We
must show that \m{u=0}. Let \m{v\in\ke_x} over $u$. Then we have
\m{tv\in\overline{\kf}_x}. Let $U$ be a neighborhood of $x$ such that \m{tv}
comes from \m{w\in\overline{\kf}(U)}. This means that \m{w_{\mid
U\backslash\pi^{-1}(P)}\in\kf(U\backslash\pi^{-1}(P))}. Since $t$ is
invertible on \m{U\backslash\pi^{-1}(P)} we can write \m{w=tw'}, with
\m{w'\in\kf(U\backslash\pi^{-1}(P))}. We have then \m{w'=v} on
\m{U\backslash\pi^{-1}(P)}. Hence \m{v\in\overline{\kf}_x} and \m{u=0}.
\end{proof}

\end{sub}

\sepsub

\Ssect{Primitive multiple curves}{cmpr}

(cf. \cite{ba_fo}, \cite{ba_ei}, \cite{dr2}, \cite{dr1}, \cite{dr4}, \cite{dr5}, 
\cite{dr6}, \cite{ei_gr}).

Let $C$ be a smooth connected projective curve. A {\em multiple curve with
support $C$} is a Cohen-Macaulay scheme $Y$ such that \m{Y_{red}=C}.

Let $n$ be the smallest integer such that \m{Y=C^{(n-1)}}, \m{C^{(k-1)}}
being the $k$-th infinitesimal neighborhood of $C$, i.e. \
\m{\ki_{C^{(k-1)}}=\ki_C^{k}} . We have a filtration \ \m{C=C_1\subset
C_2\subset\cdots\subset C_{n}=Y} \ where $C_i$ is the biggest Cohen-Macaulay
subscheme contained in \m{Y\cap C^{(i-1)}}. We call $n$ the {\em multiplicity}
of $Y$.

We say that $Y$ is {\em primitive} if, for every closed point $x$ of $C$,
there exists a smooth surface $S$, containing a neighborhood of $x$ in $Y$ as
a locally closed subvariety. In this case, \m{L=\ki_C/\ki_{C_2}} is a line
bundle on $C$ and we have \ \m{\ki_{C_j}=\ki_X^j},
\m{\ki_{C_{j}}/\ki_{C_{j+1}}=L^j} \ for \m{1\leq j<n}. We call $L$ the line
bundle on $C$ {\em associated} to $Y$. Let \m{P\in C}. Then there exist
elements $y$, $t$ of \m{m_{S,P}} (the maximal ideal of \m{\ko_{S,P}}) whose
images in \m{m_{S,P}/m_{S,P}^2} form a basis, and such that for \m{1\leq i<n}
we have \ \m{\ki_{C_i,P}=(y^{i})} .

The simplest case is when $Y$ is contained in a smooth surface $S$. Suppose
that $Y$ has multiplicity $n$. Let \m{P\in C} and \m{f\in\ko_{S,P}}  a local
equation of $C$. Then we have \ \m{\ki_{C_i,P}=(f^{i})} \ for \m{1<j\leq n},
in particular \m{I_{Y,P}=(f^n)}, and \ \m{L=\ko_C(-C)} .

We will write \m{\ko_n=\ko_{C_n}} and we will see \m{\ko_i} as a coherent sheaf
on \m{C_n} with schematic support \m{C_i} if \m{1\leq i<n}.

If $\ke$ is a coherent sheaf on $Y$ one defines its {\em generalized rank
\m{R(\ke)}} and {\em generalized degree \m{\Deg(\ke)}} (cf. \cite{dr4}, 3-): 
take any filtration of $\ke$
\[0=\ke_0\subset\ke_1\subset\cdots\subset\ke_n=\ke\]
by subsheaves such that \m{\ke_i/\ke_{i-1}} is concentrated on $C$ for 
\m{1\leq i\leq n}, then
\[R(\ke)=\sigg_{i=1}^n\text{rk}(\ke_i/\ke_{i-1})  \quad\text{and}\quad
\Deg(\ke)=\sigg_{i=1}^n\deg(\ke_i/\ke_{i-1}) .\]
Let \m{\ko_Y(1)} be a very ample line bundle on $Y$. Then the Hilbert
polynomial of $\ke$ is
\[P_\ke(m) \ = \ R(\ke)\deg(\ko_C(1))m+\Deg(\ke)+R(\ke)(1-g)\]
(where $g$ is the genus of $C$).

We deduce from proposition \ref{plongt1}:

\sepprop

\begin{subsub}\label{cmpr1}{\bf Proposition: } Let $Y$ be a multiple curve
with support $C$. Then $Y$ is a primitive multiple curve if and only if
\m{\ki_C/\ki_C^2} is zero, or a line bundle on $C$.
\end{subsub}

\sepsubsub

\begin{subsub}\label{dbl} Parametrization of double curves - \rm
In the case of double curves, D.~Bayer and D.~Eisenbud have obtained in
\cite{ba_ei} the following classification: if $Y$ is of multiplicity 2, we
have an exact sequence of vector bundles on $C$
\[0\lra L\lra\Omega_{Y\mid C}\lra\omega_C\lra 0\]
which is split if and only if $Y$ is the {\em trivial curve}, i.e. the second 
infinitesimal neighborhood of $C$, embedded by the zero section in the dual 
bundle $L^*$, seen as a surface. If $Y$ is not trivial, it is completely 
determined by the line of \m{\Ext^1_{\ko_C}(\omega_C,L)} induced by the 
preceding exact sequence. The non trivial primitive curves of multiplicity 2 
and of associated line bundle $L$ are therefore parametrized by the projective 
space
\m{\P(\Ext^1_{\ko_C}(\omega_C,L))}.
\end{subsub}

\end{sub}

\sepsub

\Ssect{Simple primitive multiple curves}{cmprq}

Let $C$ be a smooth projective irreducible curve, \m{n\geq 2} an integer and
\m{C_n} a primitive multiple curve of multiplicity $n$ and associated reduced
curve $C$. Then the ideal sheaf $\ki_C$ of $C$ in \m{C_n} is a line bundle on
\m{C_{n-1}}.

\sepprop

We say that \m{C_n} is {\em simple} if \ \m{\ki_C\simeq\ko_{n-1}}.

\sepprop

In this case the line bundle on $C$ associated to \m{C_n} is \m{\ko_C}. The
following result is proved in \cite{dr6} (th\'eor\`eme 1.2.1):

\sepprop

\begin{subsub}\label{theo_121}{\bf Theorem: }
Suppose that \m{C_n} is simple. Then there exists a flat family of smooth
projective curves \ \m{\tau:\kc\to \C} \ such that \m{\tau^{-1}(0)\simeq C}
and that \m{C_n} is isomorphic to the $n$-th infinitesimal neighborhood of
$C$ in $\kc$.
\end{subsub}

\end{sub}

\sepsec

\section{Reducible reduced deformations of primitive multiples
curves}\label{defcmp}

\Ssect{Connected Components}{compconn}

Let \m{(S,P)} be the germ of a smooth curve and \m{t\in\ko_{S,P}} a generator 
of the maximal ideal. Let \m{n>0} be an integer and \m{Y=C_n} a projective
primitive multiple curve of multiplicity $n$.

\sepprop

Let \m{k>0} be an integer. Let \ \m{\pi:\kc\to S} be a flat morphism, where
$\kc$ is a reduced algebraic variety, such that
\begin{enumerate}
\item[--] For every closed point $s\in S$ such that $s\not=P$, the fiber
$\kc_s$ has $k$ irreducible components, which are smooth and transverse, and
any three of these components have no common point.
\item[--] The fiber $\kc_P$ is isomorphic to $C_n$.
\end{enumerate}
It is easy to see that the irreducible components of $\kc$ are reduced
surfaces.

Let $Z$ be the open subset of \m{\kc\backslash\kc_P} of points $z$ belonging to
only one irreducible component of $\kc_{\pi(z)}$. Then the restriction of
\m{\pi : Z\to S\backslash\lbrace P\rbrace} is a smooth morphism. For every
\m{s\in S\backslash\lbrace P\rbrace}, let \ \m{\kc'_s=\kc_s\cap Z}. It is the
open subset of smooth points of \m{\kc_s}.

Let \m{z\in Z} and \m{s=\pi(z)}. There exist a neighborhood (for the Euclidean
topology) $U$ of $s$, isomorphic to $\C$, and a neighborhood $V$ of $z$
such that \m{V\simeq\C^2}, \m{\pi(V)=U}, the restriction of \m{\pi:V\to U}
being the projection \m{\C^2\to\C} on the first factor. We deduce easily from
that the following facts:
\begin{enumerate}
\item[--] let $s\in S\backslash\lbrace P\rbrace$ and $C_1$ an irreducible
component of \m{\kc_s}. Let $z_1,z_2\in C_1\cap Z$. Then there exist
neighborhoods (in $Z$, for the Euclidean topology) $U_1$, $U_2$ of $z_1$, $z_2$
respectively, such that if $y_1\in U_1$, $y_2\in U_2$ are such that
$\pi(y_1)=\pi(y_2)$, then $y_1$ and $y_2$ belong to the same irreducible
component of $\kc_{\pi(y_1)}$.
\item[--] for every continuous map \ $\sigma:\lbrack 0,1\rbrack\to 
S\backslash\{P\}$ \ and every $z\in Z$ such that \m{\sigma(0)=\pi(z)} there 
exists a lifting of $\sigma$, $\sigma':\lbrack 0,1\rbrack\to Z$ \ such that 
$\sigma'(0)=z$. Moreover, if \ \m{\sigma'':\lbrack 0,1\rbrack\to Z} \ is 
another lifting of $\sigma$ such that $\sigma''(0)=z$, then $\sigma'(1)$ and 
$\sigma''(1)$ are in the same irreducible component of $\kc_{\sigma(1)}$. More 
generally, if we only impose that $\sigma''(0)$ is in the same irreducible 
component of \m{\kc_{\sigma(0)}} as $z$, then $\sigma'(1)$ and $\sigma''(1)$ 
are in the same irreducible component of $\kc_{\sigma(1)}$.
\end{enumerate}

\sepprop

\begin{subsub}\label{compconn1}{\bf Lemma: } Let
\ \m{\sigma_0,\sigma_1:\lbrack 0,1\rbrack\to S\backslash\lbrace P\rbrace} \
be two continuous maps such that \Nligne \m{\sigma_0(0)=\sigma_1(0)},
\m{s=\sigma_0(1)=\sigma_1(1)}. Suppose that they are homotopic. Let
\m{\sigma'_0}, \m{\sigma'_1} be liftings \m{\lbrack 0,1\rbrack\to Z} of
\m{\sigma_0}, \m{\sigma_1} respectively, such that
\m{\sigma'_0(0)=\sigma'_1(0)}. Then \m{\sigma'_0(1)} and \m{\sigma'_1(1)}
belong to the same irreducible component of \m{\kc'_s}.
\end{subsub}

\begin{proof}
Let
\[\Psi:\lbrack 0,1\rbrack\times\lbrack 0,1\rbrack\lra S\backslash\lbrace
P\rbrace\]
be an homotopy:
\[\Psi(0,t)=\sigma_0(t),\quad \Psi(1,t)=\sigma_1(t),\quad
\Psi(t,0)=\sigma_0(0),\quad \Psi(t,1)=\sigma_0(1)\]
for \m{0\leq t\leq 1}. For every \m{u\in\lbrack 0,1\rbrack} and \m{\epsilon>0}
let \m{I_{u,\epsilon}=\lbrack u-\epsilon,u+\epsilon\rbrack\cap\lbrack
0,1\rbrack}. By using the local structure of \m{\pi_{\mid Z}} for the
Euclidean topology it is easy to see that for every \m{u\in\lbrack 0,1\rbrack},
there exists an \m{\epsilon>0} such that the restriction of $\Psi$
\[I_{u,\epsilon}\times\lbrack 0,1\rbrack\lra S\backslash\lbrace P\rbrace\]
can be lifted to a morphism
\[\Psi': I_{u,\epsilon}\times\lbrack 0,1\rbrack\lra Z\]
such that \ \m{\Psi'(t,0)=\sigma'_0(0)} \ for every \m{t\in I_{u,\epsilon}}. It
follows that if \m{I_{u,\epsilon}=\lbrack
a_{u,\epsilon},b_{u,\epsilon}\rbrack}, then \m{\Psi'(a_{u,\epsilon},1)} and
\m{\Psi'(b_{u,\epsilon},1)} are in the same irreducible component of
\m{\kc'_{\sigma_0(1)}}. We have just to cover \m{\lbrack 0,1\rbrack} with
a finite number of intervals \m{I_{u,\epsilon}} to obtain the result.
\end{proof}

\sepprop

Let \m{s\in S\backslash\lbrace P\rbrace}, \m{D_1,\ldots,D_k} be the irreducible
components of \m{\kc'_s} and \m{x_i\in D_i} for \m{1\leq i\leq k}. Let
\m{\sigma} be a loop of \m{S\backslash\lbrace P\rbrace} with origin $s$,
defining a generator of \m{\pi_1(S\backslash\lbrace P\rbrace)}. Let $i$ be an
integer such that \m{1\leq i\leq k}. The liftings \m{\sigma':\lbrack
0,1\rbrack\to Z} of $\sigma$ such that \m{\sigma'(0)=x_i} end up at a
component \m{D_j} which does not depend on \m{x_i}. Hence we can write
\[j \ = \alpha_\kc(i) .\]

\sepprop

\begin{subsub}\label{compconn2}{\bf Lemma: } \m{\alpha_\kc} is a permutation
of \m{\lbrace 1,\ldots,k\rbrace}.
\end{subsub}

\begin{proof} Suppose that \m{i\not=j} and \m{\alpha_\kc(i)=\alpha_\kc(j)}. By
inverting the paths we find liftings of paths from \m{D_{\alpha_\kc(i)}} to
\m{D_i} and \m{D_j}. This contradicts lemma \ref{compconn1}.
\end{proof}

\sepprop

Let \m{p>0} be an integer such that \m{\alpha_\kc^p=I_{\lbrace
1,\cdots,k\rbrace}}.
Let $t$ be a generator of the maximal ideal of \m{\ko_{S,P}}, $K$ the field of
rational functions on $S$ and \m{K'=K(t^{1/p})}. Let \m{S'} be the germ of the
curve corresponding to \m{K'}, \m{\theta:S'\to S} canonical the morphism and
\m{P'} the unique point of \m{\theta^{-1}(P)}. Let \m{\kd=\theta^*(\kc)}. We
have therefore a cartesian diagram
\xmat{\kd\ar[d]^\Theta\ar[r]^\rho & S'\ar[d]^\theta \\ \kc\ar[r]^\pi & S}
where $\rho$ is flat, and for every \m{s'\in S'}, $\Theta$ induces an
isomorphism \m{\kd_{s'}\simeq\kc_{\theta(s')}}. We have
\[\alpha_{\kd} \ = \ I_{\lbrace 1,\ldots,k\rbrace} \ .\]
Let \m{Z'\subset\kd} be the complement of the union of \m{\rho^{-1}(P')}
and of the singular points of the curves \m{\kd_{s'}}, \m{s'\not=P'} (hence
\m{Z'=\Theta^{-1}(Z)}).

\sepprop

\begin{subsub}\label{compconn3}{\bf Proposition: } The open subset \m{Z'} has
exactly $k$ irreducible components \m{Z'_1,\ldots,Z'_k}. Let
\m{\overline{Z'_1},\ldots,\overline{Z'_k}} be their closures in $\kd$.
Then for every \m{s'\in S'\backslash\lbrace P'\rbrace}, the
\m{Z'_i\cap\kd_{s'}}, \m{1\leq i\leq k}, are the irreducible components
of \m{\kd_{s'}} minus the intersection points with the other components, and
the \m{\overline{Z'_i}\cap\kd_{s'}} are the irreducible components of
\m{\kd_{s'}}.
\end{subsub}

\sepprop

\begin{subsub}\label{defdefred}{\bf Definition: }
Let \m{k>0} be an integer. A {\em reducible deformation of length $k$ of
\m{C_n}} is a flat morphism \ \m{\pi:\kc\to S}, where $\kc$ is a reduced
algebraic variety, such that
\begin{enumerate}
\item[--] For every closed point $s\in S$, $s\not=P$, the fiber $\kc_s$ has
$k$ irreducible components, which are smooth and transverse, and any three
of these components have no common point.
\item[--] The fiber $\kc_P$ is isomorphic to $C_n$.
\item[--] We have \ $\alpha_{\kc}=I_{\lbrace 1,\ldots,k\rbrace}$.
\end{enumerate}
\end{subsub}

\end{sub}

\sepsub

\Ssect{Maximal reducible deformations}{defmax}

Let \m{(S,P)} be the germ of a smooth curve and \m{t\in\ko_{S,P}} a generator 
of the maximal ideal. Let \m{n>0} be an integer and \m{Y=C_n} a projective
primitive multiple curve of multiplicity $n$, with underlying smooth curve
$C$. We denote by $g$ the genus of $C$ and $L$ the line bundle on $C$ 
associated to \m{C_n}.

Let \m{\pi:\kc\to S} be a reducible deformation of length $k$ of \m{C_n}. Let
\m{Z_1,\ldots,Z_k} be the closed subvarieties of
\m{\pi^{-1}(S\backslash\lbrace P\rbrace)} such that for every \m{s\in
S\backslash\lbrace P\rbrace}, \m{Z_{1s},\ldots,Z_{ks}} are the irreducible
components of \m{\kc_s} (cf. prop. \ref{compconn3}).

For \m{1\leq i\leq k}, we denote by \m{\kj_i} the ideal sheaf of \
\m{Z_1\cup\cdots\cup Z_i} \ in \m{\pi^{-1}(S\backslash\lbrace P\rbrace)}.
This sheaf is flat on  \m{S\backslash\lbrace P\rbrace}, and we have
\[0=\kj_k\subset\kj_{k-1}\subset\cdots\subset\kj_1\subset\ko_{\pi^{-1}
(S\backslash\lbrace P\rbrace)} \ .\]
The quotients \m{\ko_{\pi^{-1}(S\backslash\lbrace P\rbrace)}/\kj_1}, \m{\kj_i/
\kj_{i+1}}, \m{1\leq i<k}, are also flat on \m{S\backslash\lbrace
P\rbrace}. We obtain the filtration of sheaves on $\kc$
\[0=\overline{\kj_k}\subset\overline{\kj_{k-1}}\subset\cdots\subset
\overline{\kj_1}\subset\ko_\kc \ .\]
(cf. \ref{cohpla}). According to proposition \ref{defcmp1} the quotients
\m{\ko_\kc/\overline{\kj_1}} and \m{\overline{\kj_i}/\overline{\kj_{i+1}}},
\m{1\leq i<n}, are flat on $S$. We have \ \m{\ko_{\pi^{-1}(S\backslash\lbrace
P\rbrace)}/\kj_1 = \ko_{Z_1}}. We denote by \m{{\bf X}_i} the closed subvariety
of $\kc$ corresponding to the ideal sheaf \m{\overline{\kj_i}}.

Similarly we consider the ideal sheaf \m{\kj'_i} of \
\m{Z_{i+1}\cup\cdots\cup Z_n} \ on \m{\pi^{-1}(S\backslash\lbrace P\rbrace)},
the associated ideal sheaf \m{\overline{\kj'_i}} on $\kc$ and the
corresponding subvariety \m{{\bf X}'_i}.

\sepprop

\begin{subsub}\label{defmax1}{\bf Proposition: } We have \ \m{k\leq n} .
\end{subsub}

\begin{proof}
Let \m{\ke_0=\ko_\kc/\overline{\kj_1}} and
\m{\ke_i=\overline{\kj_i}/\overline{\kj_{i+1}}} for \m{1\leq i<n}. The
sheaves \m{\ke_{iP}} are not concentrated on a finite number of points.
To see this we use a very ample line bundle \m{\ko(1)} on $\kc$. The Hilbert
polynomial of \m{\ke_{iP}} is the same as that of \m{\ke_{is}}, \m{s\not=P},
hence it is not constant. So we have \m{R(\ke_i)\geq 1} ( \m{R(\ke_i)} is the 
generalized rank of \m{\ke_i}, cf. \ref{cmpr}), and
since
\begin{equation}\label{equ1}n \ = \ R(\ko_{C_n}) \ = \sigg_{i=0}^kR(\ke_{iP}) \
,\end{equation}
we have \m{k\leq n}.
\end{proof}

\sepprop

\begin{subsub}\label{defmax2}{\bf Definition: } We say that $\pi$ (or $\kc$)
is a {\em maximal reducible deformation} of \m{C_n} if \m{k=n}.
\end{subsub}

\sepprop

\begin{subsub}\label{defmax3}{\bf Theorem: } Suppose that $\kc$ is a maximal
reducible deformation of \m{C_n}. Then we have, for \m{1\leq i<n}
\[\overline{\kj_{i},}_P \ = \ \ki_{C_i,C_n}\]
and \m{{\bf X}_i} is a maximal reducible deformation of \m{C_i}.
\end{subsub}

\begin{proof}
Let \m{\ko_\kc(1)} be a very ample line bundle on $\kc$.

Let $Q$ be a closed point of $C$. Let \m{z\in\ko_{n,Q}} be an equation of $C$
and \m{x\in\ko_{n,Q}} over a generator of the maximal ideal of $Q$ in
\m{\ko_{C,Q}}. Let \m{\mathbf{z},\mathbf{x}\in\ko_{\kc,Q}} be over \m{z,x}
respectively. The maximal ideal of \m{\ko_{n,Q}} is \m{(x,z)}. The maximal
ideal of \m{\ko_{\kc,Q}} is generated by \m{\mathbf{z},\mathbf{x},t}. It 
follows from proposition \ref{plongt1} that there exist a neighborhood $U$ 
of $Q$ in $\kc$ and an embedding \ \m{j:U\to\P_3}. We can assume that the 
restriction of $j$ to \m{\ov{Z_1}\cap U} is 
induced by the morphism \m{\phi:\C[X,Z,T]\to\ko_{\overline{Z_1},Q}} \ of 
$\C$-algebras which associates $x$, $z$, $t$ to $X$, $Z$, $T$ respectively.

Since $\kc$ is reduced, $U$ is an open subset of a reduced hypersurface of 
$\P_3$ having $n$ irreducible components, corresponding to 
\m{\ov{Z_1},\ldots,\ov{Z_n}}. It is then clear that ${\bf X}_i$, being the 
smallest subscheme of $\kc$ containing \m{Z_1\backslash C,\ldots,Z_i\backslash 
C}, is the union in $U$ of the first $i$ hypersurface components. 

Since \m{j(\ov{Z_1})} is a hypersurface, the kernel of $\phi$ is a principal 
ideal generated by the equation $F$ of the image of $Z_1$. 

Recall that \m{\ko_n=\ko_{C_n}=(\ko_\kc)_P}.
We have \ \m{R(\ko_n/\overline{\kj_{1},}_P)=1} \ according to \eqref{equ1}.
Hence there exists a nonempty open subset $V$ of \m{C_n} such that
\m{\big(\ko_n/\overline{\kj_{1},}_P\big)_{\mid V}} is a line bundle on \m{V\cap
C}. It follows that the projection \m{\ko_n\to\ko_C} vanishes on
\m{\ov{\kj_{1},}_{P\mid V}}. Since \m{\ko_C} is torsion free this projection
vanishes everywhere on \m{\ov{\kj_1}}, i.e. \m{\ov{\kj_1}_P\subset\ki_{C,C_n}},
with equality on $V$.

The sheaf \m{\ke_0=\ko_\kc/\overline{\kj_1}} is the structural sheaf of
\m{\overline{Z_1}}, and the projection \m{\overline{Z_1}\to S} is a flat
morphism. For every \m{s\in S\backslash\lbrace P\rbrace},
\m{(\overline{Z_1})_s} is a smooth curve. The fiber \m{(\overline{Z_1})_P}
consists of $C$ and a finite number of embedded points. There exist
flat families of curves whose general fiber is smooth and the special fiber
consists of an integral curve and some embedded points (cf. \cite{ha}, III,
Example 9.8.4). We will show that this cannot happen in our case, i.e. we have
\m{\ov{\kj_1}_P=\ki_{C,C_n}}.

Let \m{\mathbf{m}=(X,Z,T)\subset\C[X,Z,T]}, and \m{\mathbf{m}_{Z_1}} the 
maximal ideal of \m{\ko_{\overline{Z_1},Q}}. 
The
ideal of \m{(\overline{Z_1})_P} in \m{\ko_{n,Q}} contains \m{z^q} and \m{x^pz}
(for suitable minimal integers \m{p\geq 0, q>0}), with \m{p>0} if and only if
$Q$ is an embedded point. Hence the ideal of \m{\overline{Z_1}} in
\m{\ko_{\kc,Q}} contains elements of type \m{\mathbf{x}^p\mathbf{z}-t\alpha},
\m{\mathbf{z}^q-t\beta}, with \m{\alpha,\beta\in \ko_{\kc,Q}}.

Let \m{\widehat{\ko_{\overline{Z_1},Q}}} be the completion of 
\m{\ko_{\overline{Z_1},Q}} with respect to \m{\mathbf{m}_{Z_1}} and
\[\widehat{\phi}:\C((X,Z,T))\lra\widehat{\ko_{\overline{Z_1},Q}}\]
the morphism deduced from $\phi$. We can also see
\m{\widehat{\ko_{\overline{Z_1},Q}}} as the completion with respect to
\m{(X,Z,T)} of \m{\ko_{\overline{Z_1},Q}} seen as a \m{\C[X,Z,T]}-module. It
follows that \ \m{\ker(\widehat{\phi})=(F)} (cf. \cite{eis}, lemma 7.15). Note
that \m{\widehat{\phi}} is surjective (this is why we use completions). Let
\m{\boldsymbol{\alpha},\boldsymbol{\beta}\in\C((X,Y,Z))} be such that
\m{\widehat{\phi}(\boldsymbol{\alpha})=\alpha},
\m{\widehat{\phi}(\boldsymbol{\beta})=\beta}. So we have
\[X^pZ-T\boldsymbol{\alpha}, \ Z^q-T\boldsymbol{\beta} \ \in \
\ker(\widehat{\phi}) \ .\]
Hence there exist \m{A,B\in\C((X,Z,T))} such that \
\m{X^pZ-T\boldsymbol{\alpha}=AF}, \m{Z^q-T\boldsymbol{\beta}=BF}. We can write
in an unique way
\[A=A_0+TA_1 , \ B=B_0+TB_1 , \ F=F_0+TF_1 ,\]
with \m{A_0,B_0,F_0\in\C((X,Z))} and \m{A_1,B_1,F_1\in\C((X,Z,T))}, and we have
\[A_0F_0=X^pZ , \ B_0F_0=Z^q \ .\]
Since $F$ is not invertible, it follows that $F_0$ is of the form \
\m{F_0=cZ}, with \m{c\in\C((X,Z,T))} invertible. So we have \m{F=cZ+TF_1}.
It follows that \m{z\in (t)} in \m{\widehat{\ko_{\overline{Z_1},Q}}}. This
implies that this is also true in \m{\ko_{\overline{Z_1},Q}}: in fact
the assertion in \m{\widehat{\ko_{\overline{Z_1},Q}}} implies that
\[z \ \in \ \bigcap_{n\geq 0}((t)+\mathbf{m}_{Z_1}^n)\]
in \m{\ko_{\overline{Z_1},Q}}, and the latter is equal to \m{(t)} according to
\cite{sa-za}, vol. II, chap. VIII, theorem 9. Hence \m{z\in(t)} in
\m{\ko_{\overline{Z_1},Q}}, i.e. \m{p=0} and  $Q$ is not an embedded point. So
there are no embedded points. This implies that
\m{\overline{\kj_{1}}_P=\ki_{C,C_n}}. Similarly, if \m{I_j} denotes the ideal 
sheaf of \m{\ov{Z_j}} for \m{1\leq j\leq n}, we have \m{I_{j,P}=\ki_{C,C_n}}.
Since the restriction of \m{\pi:\ov{Z_j}\to S} is flat, the curves 
\m{\ke_{j,s}}, \m{s\not=P}, have the same genus as $C$, and the same Hilbert 
polynomial with respect to $\ko_\kc(1)$.

Now we show that 
\m{{\bf X}'_1} is a maximal reducible deformation of \m{C_{n-1}}. We need only 
to show that \m{{\bf X}'_{1,P}=C_{n-1}}. As we have seen, for \m{2\leq j\leq 
n}, a local equation of \m{\ov{Z_j}} at any point \m{Q\in C} induces a 
generator $u_j$ of \m{\ki_{C,C_n,Q}}. Hence \m{u=\prod_{2\leq j\leq n}u_j} is 
a generator of \m{\ki_{C_{n-1},C_n,Q}}. But $u=0$ on \m{{\bf X}'_1}. It 
follows that \m{{\bf X}'_{1,P}\subset C_{n-1}}. But the Hilbert polynomial of 
\m{\ko_{C_{n-1}}} is the same as that of the structural sheaves of the fibers 
of the flat morphism \m{{\bf X}'_{1}\to S} over \m{s\not=P}, hence the same as 
\m{\ko_{{\bf X}'_{1,P}}}. Hence \m{{\bf X}'_{1,P}=C_{n-1}}.

The theorem \ref{defmax3} is then easily proved by induction on $n$.
\end{proof}

\sepprop

\begin{subsub}\label{defmax4}{\bf Corollary: } Let \m{s\in S\backslash\{P\}}
and \m{D_1,D_2} be two irreducible components of \m{\kc_s}. Then \m{D_1} is of
genus $g$ and \m{D_1\cap D_2} consists of \m{-\deg(L)} points.
\end{subsub}

\begin{proof} According to theorem \ref{defmax3}, there exists a flat family
of smooth curves $\bf C$ parametrized by $S$ such that \m{{\bf C}_P=C} and
\m{{\bf C}_s=D_1}. So the genus of \m{D_1} is equal to that of $C$.

Let us prove the second assertion. Again according to theorem \ref{defmax3} we
can suppose that \m{n=2}. We have then \ \m{\chi(\kc_s)=
\chi(C_2)=2\chi(C)+\deg(L)}. Let \m{x_1,\ldots,x_N} be the intersection points
of \m{D_1} and \m{D_2}. We have an exact sequence
\[0\lra\ko_{D_2}(-x_1-\cdots-x_N)\lra\ko_{\kc_s}\lra\ko_{D_1}\lra 0 .\]
Whence \ \m{\chi(\ko_{\kc_s})=\chi(D_1)+\chi(D_2)-N=2\chi(\ko_C)-N} (according
to the first assertion). Whence \ \m{N=-\deg(L)}.
\end{proof}

\sepprop

\begin{subsub}\label{rem0}\rm It follows from the previous results that if
\m{\pi:\kc\to S} is a maximal reducible deformation of \m{C_n}, then we have
\begin{enumerate}
\item[(i)] $\deg(L)\leq 0$ .
\item[(ii)] $\kc$ has exactly $n$ irreducible components
$\kc_1\ldots,\kc_n$.
\item[(iii)] For $1\leq i\leq n$, the restriction of $\pi$, $\pi_i:\kc_i\to S$
is a flat morphism , and \ $\pi_i^{-1}(P)=C$.
\item[(iv)] For every nonempty subset $I\subset\{1,\ldots,n\}$, let $\kc_I$ be
the union of the $\kc_i$ such that $i\in I$, and $m$ the number of elements of
$I$. Then the restriction of $\pi$, $\pi_I:\kc_I\to S$ is a maximal reducible
deformation of $C_m$.
\end{enumerate}
\end{subsub}
\sepprop

The following is immediate, and shows that we need only to consider maximal
reducible deformations parametrized by a neighborhood of 0 in $\C$:

\sepprop

\begin{subsub}\label{equ00}{\bf Proposition: } Let \m{t\in\ko_S(P)} be a
generator of the maximal ideal, and \m{\pi:\kc\to S} a maximal reducible
deformation of \m{C_n}. Let \m{S'\subset S} be an open neighborhood of $P$
where $t$ is defined and \m{\kc'=\pi^{-1}(U)}, \m{V=t(U)}. Then
\m{\pi'=t\circ\pi:\kc'\to V} is a maximal reducible deformation of \m{C_n}.
\end{subsub}

\end{sub}

\sepsec

\section{Fragmented deformations of primitive multiple curves}\label{Defecl}

The fragmented deformations of primitive multiple curves are particular cases
of reducible deformations.

In this chapter \m{(S,P)} denotes the germ of a smooth curve. Let 
\m{t\in\ko_{S,P}} be a generator of the maximal ideal of $P$. We can suppose
that $t$ is defined on the whole of $S$, and that the ideal sheaf of $P$ in
$S$ is generated by $t$.
\sepsub

\Ssect{Fragmented deformations and gluing}{eclelm}

Let \m{n>0} be an integer and \m{Y=C_n} a projective primitive multiple curve 
of multiplicity $n$.

\sepprop

\begin{subsub}\label{defdefecl}{\bf Definition: } Let \m{k>0} be an integer.
A {\em general fragmented deformation of length $k$} of \m{C_n} is a flat
morphism \ \m{\pi:\kc\to S} \ such that for every point \m{s\not=P} of $S$,
the fiber \m{\kc_s} is a disjoint union of $k$ projective smooth irreducible
curves, and such that \m{\kc_P} is isomorphic to \m{C_n}.
\end{subsub}

\sepprop

We have then \m{k\leq n}. If \m{k=n} we say that $\pi$ (or $\kc$) is a {\em
general maximal fragmented deformation} of \m{C_n}. We suppose in the sequel
that it is the case.

The line bundle on $C$ associated to \m{C_n} is \m{\ko_C} (by proposition
\ref{defmax4}).

Let \m{p>0} be an integer. Let \m{K} be the field of rational functions
on $S$  and \m{K'=K(t^{1/p})}. Let \m{S'} be the germ of curve corresponding
to \m{K'}, \m{\theta:S'\to S} the canonical morphism and \m{P'} the unique
point of \m{\theta^{-1}(P)}. Let \m{\kd=\theta^*(\kc)}. So we have a cartesian
diagram
\xmat{\kd\ar[d]^\Theta\ar[r]^\rho & S'\ar[d]^\theta \\ \kc\ar[r]^\pi & S}
where $\rho$ is flat, and for every \m{s'\in S'}, $\Theta$ induces an
isomorphism \m{\kd_{s'}\simeq\kc_{\theta(s')}}.

\sepprop

\begin{subsub}\label{ecl1}{\bf Proposition: } For a suitable choice of $p$,
$\kd$ has exactly $n$ irreducible components \m{\kd_1,\ldots,\kd_n}, and for
every point \m{s\not=P'} of \m{S'}, \m{\kd_{1s},\ldots,\kd_{ns}} are the 
irreducible components of \m{\kd_s}, for \m{1\leq i\leq n} the restriction
of $\rho$: \m{\kd_{is}\to S'} is flat, and  \m{\kd_{P'}=C_n}.
\end{subsub}

(See proposition \ref{compconn3})

\sepprop

\begin{subsub}\label{defdefecl2}{\bf Definition: } A {\em fragmented
deformation} of \m{C_n} is a general maximal fragmented deformation of length
$n$ of \m{C_n} having $n$ irreducible components.
\end{subsub}

\sepprop

We suppose in the sequel that $\kc$ is a fragmented deformation of \m{C_n},
union of $n$ irreducible components \m{\kc_1,\ldots,\kc_n}.

\sepprop

\begin{subsub}\label{ecl2}{\bf Proposition: } Let \m{I\subset\{1,\ldots,n\}} be
a nonempty subset having $m$ elements. Let \m{\kc_I=\cup_{i\in I}\kc_i}. Then
the restriction of $\pi$, \m{\kc_I\to S}, is flat, and the fiber \m{\kc_{IP}}
is canonically isomorphic to \m{C_m}.
\end{subsub}

(See \ref{rem0})

\sepprop

In particular there exists a filtration of ideal sheaves
\[0\subset\ki_1\subset\cdots\subset\ki_{n-1}\subset\ko_\kc\]
such that for \m{1\leq i<n} and \m{s\in S\backslash\{P\}}, \m{\ki_{is}} is the
ideal sheaf of \m{\cup_{j=i}^n\kc_{js}}, and that \m{\ki_{iP}} is that of
\m{C_{n-i}}.

\sepprop

\begin{subsub}\label{ecl2c}{\bf Definition: } For \m{1\leq i\leq n}, let
\m{\pi_i:\kc_i\to S} be a flat family of smooth projective irreducible curves,
with a fixed isomorphism \m{\pi_i^{-1}(P)\simeq C}. A {\em gluing of
\m{\kc_1,\cdots,\kc_n} along $C$} is an algebraic variety $\kd$ such that
\begin{enumerate}
\item[-] for $1\leq i\leq n$, $\kc_i$ is isomorphic to a closed subvariety of
$\kd$, also denoted by $\kc_i$, and $\kd$ is the union of these subvarieties.
\item[-] $\coprod_{1\leq i\leq n}(\kc_i\backslash C)$ is an open subset of
$\kd$.
\item[-] There exists a morphism \ $\pi:\kd\to S$ \ inducing $\pi_i$ on
$\kc_i$, for $1\leq i\leq n$.
\item[-] The subvarieties \ $C=\pi_i^{-1}(P)$ \ of $\kc_i$ coincide in $\kd$.
\end{enumerate}
\end{subsub}

\sepprop

For example the previous fragmented deformation $\kc$ of $C_n$ is a gluing
of \m{\kc_1,\cdots,\kc_n} along $C$.

All the gluings of \m{\kc_1,\cdots,\kc_n} along $C$ have the same underlying
Zariski topological space.

Let $\BS{\ka}$ be the {\em initial gluing} of the \m{\kc_i} along $C$. It is an
algebraic variety whose points are the same as those of $\kc$, i.e.
\[(\coprod_{i=1}^n\kc_i)/\sim \quad ,\]
where $\sim$ is the equivalence relation: if \m{x\in\kc_i} and
\m{y\in\kc_j}, \m{x\sim y} if and only if \m{x=y}, or if
\m{x\in\kc_{iP}\simeq C}, \m{y\in\kc_{jP}\simeq C} and \m{x=y} in $C$.
The structural sheaf is defined by: for every open subset $U$ of $\BS{\ka}$
\[\ko_\BS{\ka}(U) \ = \
\{(\alpha_1,\ldots,\alpha_n)\in\ko_{\kc_1}(U\cap\kc_1)\times
\cdots\times\ko_{\kc_n}(U\cap\kc_n) ; \alpha_{1\mid C}=\cdots=\alpha_{n\mid 
C}\} .\]
For every gluing $\kd$ of \m{\kc_1,\cdots,\kc_n}, we have an obvious dominant
morphism \ \m{\BS{\ka}\to\kd}. If follows that the sheaf of rings \m{\ko_\kd}
can be seen as a subsheaf of \m{\ko_{\BS{\ka}}}.

The fiber \m{D=\BS{\ka}_P} is not a primitive multiple curve (if \m{n>2}): if
\m{\ki_{C,D}} denotes the ideal sheaf of $C$ in $D$ we have \m{\ki_{C,D}^2=0},
and \m{\ki_{C,D}\simeq\ko_C\ot\C^{n-1}} .

\sepprop

\begin{subsub}\label{ecl2b}{\bf Proposition: } Let $\kd$ be a gluing of
\m{\kc_1,\cdots,\kc_n}. Then
\m{\pi^{-1}(P)} is a primitive multiple curve if and only if for every closed
point $x$ of $C$, there exists a neighborhood of $x$ in $\kd$ that can be
embedded in a smooth variety of dimension 3.
\end{subsub}

\begin{proof} Suppose that \m{\pi^{-1}(P)} is a primitive multiple curve. Then
\m{\ki_C/(\ki_C^2+(\pi))} is a principal module at $x$: suppose that the
image of \m{u\in m_{\kd,x}} is a generator. The module \m{m_{\kd,x}/\ki_C} is
also principal (since it is the maximal ideal of $x$ in $C$): suppose that
the image of \m{v\in m_{\kd,x}} is a generator. Then the images of \m{u,v,\pi}
generate \m{m_{\kd,x}/m_{\kd,x}^2}, so according to proposition \ref{plongt1},
we can locally embed $\kd$ in a smooth variety of dimension 3.

Conversely, suppose that a neighborhood of \m{x\in C} in $\kd$ is embedded in
a smooth variety $Z$ of dimension 3. The proof of the fact that 
\m{\pi^{-1}(P)} is Cohen-Macaulay is similar to that of theorem 3.2.3. We can 
suppose that $\pi$ is defined on $Z$. We have \ 
\m{\pi_{\mid\kc_1}=\pi_1\not\in m_{\kc_1,x}^2}, so \m{\pi\not\in m_{Z,x}^2}. 
It follows that the surface of $Z$ defined by $\pi$ is smooth at $x$, and that 
we can locally embed \m{\pi^{-1}(P)} in a smooth surface. Hence 
\m{\pi^{-1}(P)} is a primitive multiple curve.
\end{proof}

\end{sub}

\sepsub

\Ssect{Fragmented deformations of length 2}{ecl_l2}

Let \m{\pi:\kc\to S} be a fragmented deformation of \m{C_2}. So $\kc$ has two 
irreducible components \m{\kc_1}, \m{\kc_2}. Let $\ka$ be the initial gluing 
of \m{\kc_1} and \m{\kc_2} along $C$. For every open subset $U$ of $\kc$, $U$ 
is also an open subset of $\ka$ and \m{\ko_\kc(U)} is
a sub-algebra of \m{\ko_\ka(U)}. For \m{i=1,2}, let \m{\pi_i:\kc_i\to S} be
the restriction of $\pi$. We will also denote \m{t\circ\pi} by $\pi$, and
\m{t\circ\pi_i} by \m{\pi_i}. So we have
\m{\pi=(\pi_1,\pi_2)\in\ko_\kc(\kc)}.

Let \m{\ki_C} be the ideal sheaf of $C$ in $\kc$. Since \m{C_2=\pi^{-1}(P)} we
have \ \m{\ki_C^2\subset\span{(\pi_1,\pi_2)}} .

Let \m{m>0} be an integer, \m{x\in C}, \m{\alpha_1\in\ko_{\kc_1,x}},
\m{\alpha_2\in\ko_{\kc_2,x}}. We denote by \m{[\alpha_1]_m} (resp.
\m{[\alpha_2]_m}) the image of \m{\alpha_1} (resp. \m{\alpha_2}) in
\m{\ko_{\kc_1,x}/(\pi_1^m)} (resp. \m{\ko_{\kc_2,x}/(\pi_2^m)}).

\sepprop

\begin{subsub}\label{ecl3}{\bf Proposition:  1 -- }There exists an unique
integer \m{p>0} such that \m{\ki_C/\span{(\pi_1,\pi_2)}} is generated by the
image of \m{(\pi_1^p,0)}.

{\bf 2 -- }The image of \m{(0,\pi_2^p)} generates
\m{\ki_C/\span{(\pi_1,\pi_2)}}.

{\bf 3 -- }For every \m{x\in C}, \m{\alpha\in\ko_{\kc_1,x}} and
\m{\beta\in\ko_{\kc_2,x}}, we have \ \m{(\pi_1^p\alpha,0)\in\ko_{\kc,x}} \ and
\ \m{(0,\pi_2^p\beta)\in\ko_{\kc,x}}.
\end{subsub}

\begin{proof} Let \m{x\in C} and \m{u=(\pi_1\alpha,\pi_2\beta)} whose image is
a generator of \m{\ki_C/\span{(\pi_1,\pi_2)}} at $x$
(\m{\ki_C/\span{(\pi_1,\pi_2)}} is a locally free sheaf of rank 1 of
\m{\ko_C}-modules). Let \m{\beta_0\in\ko_{\kc_1,x}} be such that
\m{(\beta_0,\beta)\in\ko_{\kc,x}}. Then the image of
\[u-(\pi_1,\pi_2)(\beta_0,\beta)=(\pi_1(\alpha-\beta_0),0)\]
is also a generator of \m{\ki_C/\span{(\pi_1,\pi_2)}} at $x$. We can write it
\m{(\pi_1^p\lambda,0)}, where $\lambda$ is not a multiple of $\pi_1$.

Now we show that $p$ is the smallest integer $q$ such that
\m{(\ki_C/\span{(\pi_1,\pi_2)})_x} contains the image of an element of the
form \m{(\pi_1^q\mu,0)}, with $\mu$ not divisible by $\pi_1$. We can write
\[(\pi_1^q\mu,0) \ = \
(u_1,u_2)(\pi_1^p\lambda,0)+(v_1,v_2)(\pi_1,\pi_2) \]
with \m{(u_1,u_2),(v_1,v_2)\in\ko_{\kc,x}}. So we have \m{v_2=0}, hence
\m{(v_1,v_2)\in\ki_{C,x}}. So we can write \m{(v_1,v_2)} as the sum of a
multiple of \m{(\pi_1^p\lambda,0)} and a multiple of \m{(\pi_1,\pi_2)}.
Finally we obtain \m{(\pi_1^q\mu,0)} as
\[(\pi_1^q\mu,0) \ = \
(u_{12},u_{22})(\pi_1^p\lambda,0)+(v_{11},0)(\pi_1,\pi_2)^2 .\]
In the same way we see that \m{(\pi_1^q\mu,0)} can be written as
\[(\pi_1^q\mu,0) \ = \
(u_{1p},u_{2p})(\pi_1^p\lambda,0)+(v_{1p},0)(\pi_1,\pi_2)^p ,\]
which implies immediately that \m{q\geq p}.

It follows that $p$ does not depend on $x$ and that
\m{\ki_C/\span{(\pi_1,\pi_2)}} is a subsheaf of\Nligne
\m{\span{(\pi_1^p,0)}/\span{(\pi_1^{p+1},0)}\simeq\ko_C}. Since
\m{\ki_C/\span{(\pi_1,\pi_2)}} is of degree 0 by (by corollary \ref{defmax4}) 
it follows that
\m{\ki_C/\span{(\pi_1,\pi_2)}\simeq\span{(\pi_1^p,0)}/\span{(\pi_1^{p+1},0)}},
from which we deduce assertion 1- of proposition \ref{ecl3}. The second
assertion comes from the fact that \ \m{(0,\pi_2^p)=\pi^p-(\pi_1^p,0)}.

To prove the third, we use the fact that there exists
\m{\alpha'\in\ko_{\kc_2,x}} such that \m{(\alpha,\alpha')\in\ko_{\kc,x}}
(because \m{\kc_1\subset\kc}). Hence \
\m{(\pi_1^p,0)(\alpha,\alpha')=(\pi_1^p\alpha,0)\in\ko_{\kc,x}}. Similarly, we
obtain that \m{(0,\pi_2^p\beta)\in\ko_{\kc,x}}.
\end{proof}

\sepprop

According to the proof of proposition \ref{ecl3}, for every \m{x\in C}, $p$
is the smallest integer $q$ such that there exists an element of
\m{\ko_{\kc,x}} of the form \m{(\pi_1^q\alpha,0)} (resp.
\m{(0,\pi_2^q\alpha)}), with \m{\alpha\in\ko_{\kc_1,x}} (resp.
\m{\alpha\in\ko_{\kc_2,x}}) not vanishing on $C$.

Let \m{x\in C} and \m{\alpha_1\in\ko_{\kc_1,x}}. Since \m{\kc_1\subset\kc}
there  exists \m{\alpha_2\in\ko_{\kc_2,x}} such that
\m{(\alpha_1,\alpha_2)\in\ko_{\kc,x}}. Let  \m{\alpha'_2\in\ko_{\kc_2,x}} such
that \m{(\alpha_1,\alpha'_2)\in\ko_{\kc,x}}. We have then
\m{(0,\alpha_2-\alpha'_2)\in\ko_{\kc,x}}. So there exists
\m{\alpha\in\ko_{\kc_2,x}} such that \ \m{\alpha_2-\alpha'_2=\pi_2^p\alpha}. It
follows that the image of \m{\alpha_2} in \m{\ko_{\kc_2,x}/(\pi_2^p)} is
uniquely determined. Hence we have:

\sepprop

\begin{subsub}\label{ecl4}{\bf Proposition: } There exists a canonical
isomorphism
\[\Phi : \kc_1^{(p)}\lra \kc_2^{(p)}\]
between the infinitesimal neighborhoods of order \m{p} of \m{\kc_1} and
\m{\kc_2} (i.e. \m{\ko_{ \kc_i^{(p)}}=\ko_{\kc_i}/(\pi^p_i)}), such that for
every \m{x\in C}, \m{\alpha_1\in\ko_{\kc_1,x}} and
\m{\alpha_2\in\ko_{\kc_2,x}}, we have  \ \m{(\alpha_1,\alpha_2)\in\ko_{\kc,x}}
\ if and only if \ \m{\Phi_x([\alpha_1]_p)=[\alpha_2]_p}. For every \
\m{\alpha\in\ko_{\kc_1,x}} \ we have \ \m{\Phi_x(\alpha)_{\mid C}=\alpha_{\mid
C}}, and \ \m{\Phi_x(\pi_1)=\pi_2}.
\end{subsub}

\sepprop

The simplest case is \m{p=1}. In this case \m{\Phi:C\to C} is  the identity
and \m{\kc=\ka} (the {\em initial} gluing).

\sepprop

\begin{subsub}\label{recipr}Converse - \rm Recall that $\BS{\ka}$ denotes the
initial gluing of \m{\kc_1,\kc_2} (cf. \ref{ecl2c}). Let \ \m{\Phi :
\kc_1^{(p-1)}\to \kc_2^{(p-1)}} \ be an isomorphism inducing the identity on
$C$ and such that \m{\Phi(\pi_1)=\pi_2}. We define a subsheaf of algebras
\m{\ku_\Phi} of \m{\ko_\BS{\ka}}: \m{\ku_\Phi=\ko_\BS{\ka}} on
\m{\BS{\ka}\backslash C}, and for every point $x$ of $C$
\[\ku_{\Phi,x} \ = \ \{(\alpha_1,\alpha_2)\in\ko_{\kc_1,x}\times\ko_{\kc_2,x}
\ ; \ \Phi_x([\alpha_1]_p)=[\alpha_2]_p\}  \ .\]
It is easy to see that \m{\ku_\Phi} is the structural sheaf of an algebraic
variety \m{\ka_\Phi}, that the inclusion \m{\ku_\Phi\subset\ko_\BS{\ka}}
defines a dominant morphism
$\BS{\ka}\to\ka_\Phi$
inducing an isomorphism between the underlying topological spaces (for the
Zariski topology), and that the composed morphisms
\m{\kc_i\subset\BS{\ka}\to\ka_\Phi}, \m{i=1,2}, are immersions. Moreover, the
morphism \m{\pi:\BS{\ka}\to S} factorizes through \m{\ka_\Phi} :
\xmat{\BS{\ka}\ar[r]\ar@/_2pc/[rr]^\pi & \ka_\Phi\ar[r]^-{\pi_\Phi} & S}
and \ \m{\pi_\Phi:\ka_\Phi\to S} \ is flat.

For \m{2\leq i\leq p}, let \ \m{\Phi^{(i)}:\kc_1^{(i)}\to\kc_2^{(i)}} \ be
the isomorphism induced by $\Phi$.
\end{subsub}

\sepprop

\begin{subsub}\label{ecl5}{\bf Proposition: } \m{\pi_\Phi^{-1}(P)} is a
primitive double curve.
\end{subsub}

\begin{proof}
Let $x$ be a closed point of $C$. We first show that \m{\ki_{C,x}^2\subset
(\pi)}. Let \ \m{u=(\pi_1\alpha,\pi_2\beta)\in\ki_{C,x}}. Let
\m{\beta'\in\ko_{\kc_2,x}} be such that \m{\Phi_x([\alpha]_p)=[\beta']_p}. We
have then \ \m{v=(\alpha,\beta')\in\ko_{\kc,x}}. We have \ \m{u-\pi
v=(0,\pi_2(\beta-\beta'))\in\ko_{\kc,x}}. Therefore \
\m{[\pi_2(\beta-\beta')]_p=\Phi_x(0)=0}. Hence
\Nligne\m{\pi_2(\beta-\beta')\in(\pi_2^p)}. We can then write
\[u \ = \ \pi v+(0,\pi_2^p\gamma) .\]
Let \m{u'\in\ki_{C,x}}, that can be written as \ \m{u'=\pi v'+
(0,\pi_2^p\gamma')}. We have then
\[uu' \ = \ \pi.\big(\pi vv'+(0,\pi_2\gamma')v+(0,\pi_2\gamma)v'
+(0,\pi_2^{2p-1}\gamma\gamma')\big) \ \in \ (\pi) .\]
It remains to show that \ \m{\ki_{C,x}/(\pi)\simeq\ko_{C,x}}. We have
\begin{eqnarray*}
\ki_{C,x} & = & \{(\pi_1\alpha,\pi_2\beta)\in\ko_{\kc_1,x}\times\ko_{\kc_2,x};
\Phi_x([\pi_1\alpha]_p)=[\pi_2\beta]_p\}\\
& = & \{(\pi_1\alpha,\pi_2\beta)\in\ko_{\kc_1,x}\times\ko_{\kc_2,x};
\Phi_x^{(p-1)}([\alpha]_{p-1})=[\beta]_{p-1}\} , \\
(\pi)_x & = & \{(\pi_1\alpha,\pi_2\beta)\in\ko_{\kc_1,x}\times\ko_{\kc_2,x};
\Phi_x([\alpha]_p)=[\beta]_p\} .
\end{eqnarray*}
So if \ \m{(\pi_1\alpha,\pi_2\beta)\in\ki_{C,x}}, we have \
\m{w=\Phi_x([\alpha]_p)-[\beta]_p\in(\pi_2^{p-1})_x/(\pi_2^p)_x\simeq\ko_{C,x}}.
Hence we have a morphism of \m{\ko_{\kc,x}}-modules
\[\xymatrix@R=4pt{\lambda:\ki_{C,x}\ar[r] & \ko_{C,x}\\
(\pi_1\alpha,\pi_2\beta)\fmaps[r] & w}\]
whose kernel is \m{(\pi)_x}. We have now only to show that $\lambda$ is
surjective, which follows from the fact that \ \m{\lambda(\pi_1^p,0)=1}.
\end{proof}

\end{sub}

\sepsub

\Ssect{Spectrum of a fragmented deformation and ideals of
sub-deformations}{spectre}

Let \ \m{\pi:\kc\to S} \ be a fragmented deformation of \m{C_n},
\m{\kc_1,\ldots,\kc_n} the irreducible components of $\kc$. For \m{1\leq i\leq
n}, let \m{\pi_i=\pi_{\mid\kc_i}}. As in \ref{ecl_l2}, we denote also
\m{t\circ\pi_i} by \m{\pi_i}.  Let \m{I=\{i,j\}}
be a subset of \m{\{1,\ldots,n\}}, with \m{i\not=j}. Then \m{\pi:\kc_I\to S}
is a fragmented deformation of \m{C_2}. According to \ref{ecl_l2} there exists
a unique integer \m{p>0} such that \m{\ki_{C,\kc_I}/(\pi)} is generated by the
image of \m{(\pi_i^p,0)} (and also by the image of \m{(0,\pi_j^p)}). Recall
that $p$ is the smallest integer $q$ such that \m{\ki_{C,\kc_I}} contains a
non zero element of the form \m{(\pi_i^q\lambda,0)} (or \m{(0,\pi_j^q\mu)}),
with \m{\lambda_{\mid C}\not=0} (resp. \m{\mu_{\mid C}\not=0}). Let
\[p_{ij} \ = \ p_{ji} \ = \ p ,\]
and \m{p_{ii}=0} for \m{1\leq i\leq n}. The symmetric matrix
\m{(p_{ij})_{1\leq i,j\leq n}} is called the {\em spectrum} of $\kc$.

\sepprop

\begin{subsub}\label{genecl} Generators of \m{(\ki_C^p+(\pi))/
(\ki_C^{p+1}+(\pi))} - \rm Let \m{i,j\in\{1,\ldots,n\}} be such that
\m{i\not=j}. Let \m{x\in C}. Since \m{\kc_{\{i,j\}}\subset\kc} there exists an
element \ \m{{\bf u}_{ij}=(u_m)_{1\leq m\leq n}} \ of \m{\ko_{\kc,x}} such
that \m{u_i=0} and \m{u_j=\pi_j^{p_{ij}}}. According to proposition \ref{ecl2},
the image of \m{{\bf u}_{ij}} generates \m{\ki_C/(\ki_C^2+(\pi))} at $x$.
\end{subsub}

According to proposition \ref{ecl3} and the fact that the image of \m{{\bf
u}_{ij}} generates\Nligne
\m{\ki_{C,\kc_{ij},x}/(\ki_{C,\kc_{ij},x}^2+(\pi))}, for every integer
$m$ such that $m\not=i,j$ and that \m{1\leq m\leq n}, \m{u_m} is of the form
\ \m{u_m=\alpha_{ij}^{(m)}\pi_m^{p_{im}}}, with \
\m{\alpha_{ij}^{(m)}\in\ko_{\kc_m,x}} \ invertible. Let \
\m{\alpha_{ij}^{(i)}=0} \ and \ \m{\alpha_{ij}^{(j)}=1}.

\sepprop

\begin{subsub}\label{ecl7}{\bf Proposition: } {\bf 1 -- }\m{\alpha_{ij\mid
C}^{(m)}} is a non zero constant, uniquely determined and independent of $x$.

{\bf 2 -- }Let \ \m{{\bf a}_{ij}^{(m)}=\alpha_{ij\mid C}^{(m)}\in\C}.
Then we have, for all integers \m{i,j,k,m,q} such that\Nligne \m{1\leq
i,j,k,m,q\leq n}, \m{i\not=j}, \m{i\not=k}
\[{\bf a}_{ik}^{(m)}{\bf a}_{ij}^{(q)} \ = \ {\bf a}_{ik}^{(q)}{\bf
a}_{ij}^{(m)} .\]
In particular we have \ \m{{\bf a}_{ij}^{(m)}={\bf a}_{ik}^{(m)}{\bf
a}_{ij}^{(k)}} \ and \ \m{{\bf a}_{ij}^{(m)}{\bf a}_{im}^{(j)}=1}.
\end{subsub}

\begin{proof}
Let \m{{\bf u}'_{ij}} have the same properties as \m{{\bf u}_{ij}}. Then
\m{{\bf v}={\bf u}'_{ij}-{\bf u}_{ij}\in\ki_{C,x}^2+(\pi)}. So the
image of $\bf v$ in \m{\ko_{\kc_{im,}x}} belongs to \
\m{\ki_{C,\kc_{im},x}^2+(\pi)}. It follows that the $m$-th component of $\bf v$
is a multiple of \m{\pi_m^{p_{im}+1}}. Hence \m{\alpha_{ij\mid C}^{(m)}} is
uniquely determined. It follows that when $x$ varies the \m{\alpha_{ij\mid
C}^{(m)}} can be glued together and define a global section of \m{\ko_C},
which must be a constant. This proves 1-.

Now we prove 2-. There exists \m{u\in\ko_{\kc,x}} such that the $k$-th
component of $u$ is \m{\alpha^{(k)}_{ij}}, and $u$ is invertible. Then the
image of \m{(v_m)=\frac{{\bf u}_{ij}}{u}} generates \m{\ki_C/(\ki_C^2+(\pi))},
and \m{v_k=1}. Hence according to 1-, we have \ \m{v_{m\mid C}={\bf
a}_{ik}^{(m)}}, i.e.
\[\frac{{\bf a}_{ij}^{(m)}}{{\bf a}_{ij}^{(k)}} \ = \ {\bf a}_{ik}^{(m)} .\]
We have the same equality with $q$ instead of $m$, whence 2- is easily deduced.
\end{proof}

\sepprop

Let $p$ be an integer such that \m{1\leq p<n}, and 
\m{(i_1,j_1),\ldots,(i_p,j_p)}
$p$ pairs of distinct integers of \m{\{1,\ldots,n\}}. Then the image of \
\m{\prod_{m=1}^p{\bf u}_{i_mj_m}} \ is a generator of \
\m{(\ki_C^p+(\pi))/(\ki_C^{p+1}+(\pi))}.

Let \ \m{I\subset\{1,\cdots,n\}} \ be a nonempty subset, distinct from
\m{\{1,\cdots,n\}}. Let \ \m{i\in\{1,\cdots,n\}\backslash I}. Let
\[{\bf u}_{I,i} \ = \ \prod_{j\in I} {\bf u}_{ji} \ .\]
Recall that \ \m{\kc_I=\cup_{j\in I}\kc_j\subset\kc}.

\sepprop

\begin{subsub}\label{ecl8b}{\bf Proposition: } The ideal sheaf of \m{\kc_I} is
generated by \m{{\bf u}_{I,i}} at $x$.
\end{subsub}

\begin{proof} According to proposition \ref{ecl2b} there exists an embedding 
of a neighborhood of $x$ in a smooth variety of dimension 3. In this variety 
each \m{\kc_i} is a smooth surface defined by a single equation. The ideal of 
the union of the \m{\kc_i}, \m{i\in I} is the product of these equations.
\end{proof}

\sepprop

\begin{subsub}\label{ecl6}{\bf Proposition: } Let \m{i,j,k} be distinct
integers such that \m{1\leq i,j,k\leq n}. Then if \m{p_{ij}<p_{jk}}, we have
\m{p_{ik}=p_{ij}}.
\end{subsub}

\begin{proof}
We can come down to the case \m{n=3} by considering \m{\kc_{\{i,j,k\}}}. We
can suppose that \ \m{p_{23}\leq p_{12}\leq p_{13}}, and we must show that \
\m{p_{23}=p_{12}}. We have
\[{\bf u}_{21}=(\pi_1^{p_{12}},0,\alpha^{(3)}_{21}\pi_3^{p_{23}}) , \quad
{\bf u}_{31}=(\pi_1^{p_{13}},\alpha^{(2)}_{31}\pi_2^{p_{23}},0) .\]
So
\[{\bf u}_{31}-\pi^{p_{13}-p_{12}}{\bf u}_{21}=\big(0,
\alpha^{(2)}_{31}\pi_2^{p_{23}},-\alpha^{(3)}_{21}\pi_3^{p_{23}+p_{13}-p_{12}}
\big) \ \in \ \ko_{\kc,x} .\]
Taking the image of this element in \m{\ko_{\kc_{12}x}}, we see that \
\m{p_{23}\geq p_{12}}, hence \ \m{p_{23}=p_{12}}.
\end{proof}

\sepprop

\begin{subsub}\label{ecl12}{\bf Proposition: }{\bf 1 --} Let $i$, $j$ be
distinct integers such that \m{1\leq i,j\leq n}. Then we have \
\m{\ki_{C,x}=({\bf u}_{ij})+(\pi)}.

{\bf 2 --} Let \m{v=(v_m)_{1\leq m\leq n}\in\ki_{C,x}} such that \m{v_i} is a
multiple of \m{\pi_i^p}, with \m{p>0}. Then we have \ \m{v\in({\bf
u}_{ij})+(\pi^p)}.
\end{subsub}

\begin{proof} Let \ \m{N=1+\max_{1\leq k\leq n}(q_i)}, where \ 
\m{q_i=\sigg_{j=1}^np_{ij}}. For every integer $j$ 
such that \m{1\leq j\leq n} we have
\m{(0,\ldots,0,\pi_j^{q_j},0,\ldots,0)\in\ko_\kc(\kc)}. Hence \
\m{\ki_C^N\subset(\pi)}. We will show by induction on $k$ that \ \m{\ki_{C,x}
\subset({\bf u}_{ij})+(\pi)+\ki_{C,x}^k}. Taking \m{k=N} we obtain 1-.

For \m{k=1} it is obvious. Suppose that it is true for \m{k-1\geq 1}. It is
enough to prove that \ \m{\ki_{C,x}^{k-1}\subset({\bf
u}_{ij})+(\pi)+\ki_{C,x}^k}.
Let \ \m{w_1,\ldots,w_{k-1}\in\ki_{C,x}}. Since the image of \m{{\bf u}_{ij}}
generates \m{\ki_{C,x}/(\ki_{C,x}^2+(\pi))}, we can write \m{w_p} as
\[w_p \ = \ \lambda_p{\bf u}_{ij}+\pi\mu_p+\nu_p ,\]
with \m{\lambda_p,\mu_p\in\ko_{\kc,x}} and \m{\nu_p\in\ki_{C,x}^2}. So we have
\[w_1\cdots w_{k-1} \ = \ \lambda{\bf u}_{ij}+\pi\mu+\nu ,\]
with \m{\lambda,\mu\in\ko_{\kc,x}} and \m{\nu\in\ki_{C,x}^{2k-2}}. Since
\m{2k-2\geq k}, we have \ \m{w_1\cdots w_{k-1}\in({\bf
u}_{ij})+(\pi)+\ki_{C,x}^k}. This proves 1-.

We prove 2- by induction on $p$. The case $p=1$ follows 1-. Suppose that it is
true for \m{p-1\geq 1}. So we can write $v$ as
\[v \ = \ \lambda{\bf u}_{ij}+\pi^{p-1}\mu ,\]
with \m{\lambda,\mu\in\ko_{\kc,x}}. We can write \m{v_i} as 
\m{v_i=\alpha\pi^p}. So we have \
\m{\alpha\pi_i^p=\pi_i^{p-1}\mu_i}, whence \ \m{\mu_i=\alpha\pi_i}. Hence \
\m{\mu\in\ki_{C x}}. According to 1- we can write $\mu$ as \
\m{\mu=\theta{\bf u}_{ij}+\pi\tau}, with \m{\theta,\tau\in\ko_{\kc,x}}. So
\[v \ = \ (\lambda+\pi^{p-1}\theta){\bf u}_{ij}+\pi^p\tau ,\]
which proves the result for $p$.
\end{proof}

\sepprop

\begin{subsub}\label{ecl16} The ideal sheaves \m{\ki_{\kc_I}} -- \rm Recall
that \ \m{I\subset\{1,\cdots,n\}} \ is a nonempty subset, distinct from
\m{\{1,\cdots,n\}}. For every subset $J$ of \m{\{1,\cdots,n\}}, let \
\m{J^c=\{1,\cdots,n\}\backslash J} \ and \m{\ko_J=\ko_{\kc_J}}. It
follows from proposition \ref{ecl8b} that \m{\ki_{\kc_I}} is a line bundle on
\m{\kc_{I^c}}.
\end{subsub}

\sepprop

From now on, we suppose that \m{S\subset\C} and \m{P=0} (cf. proposition
3.2.6).

\sepprop

\begin{subsub}\label{ecl19}{\bf Theorem: } We have \
\m{\ki_{\kc_I}\simeq\ko_{I^c}}.
\end{subsub}

\begin{proof} By induction on $n$. If \m{n=2} the result follows from
proposition \ref{ecl3} and the fact that \m{S\subset\C}. Suppose that it is
true for \m{n-1\geq 2}. We will prove that it is true for $n$ by induction on
the number of elements $q$ of \m{I^c}. Suppose first that \m{q=1} and let $i$
be the unique element of \m{I^c}. Then according to proposition \ref{ecl8b},
\m{\ki_{\kc_I}} is generated by \m{(0,\ldots,0,\pi_i^{q_i},0,\ldots,0)}, so
the result is true in this case. Suppose that it is true if \ \m{1\leq q<k<n},
and that \m{q=k}. Let \ \m{K=\{1,\cdots,n-1\}}. We can assume that \
\m{I\subset K}.

According to proposition \ref{ecl8b}, we have, for every \m{x\in C},
\m{\ki_{\kc_I,x}\simeq\ko_{I^cx}}. We have \ \m{\ki_{\kc_K}\subset\ki_{\kc_I}},
and \ \m{\ki_{\kc_K}\simeq\ko_{\{n\}}}. We have
\[\ki_{\kc_I}/\ki_{\kc_K} \ = \ \ki_{\kc_I,\kc_K} \]
(the ideal sheaf of \m{\kc_I} in \m{\kc_K}). From the first induction
hypothesis we have
\[\ki_{\kc_I,\kc_K} \ \simeq \ \ko_{(I\cup\{n\})^c} .\]
So we have an exact sequence of sheaves
\[0\lra\ko_{\{n\}}\lra\ki_{\kc_I}\lra\ko_{I^c\backslash\{n\}}\lra 0 .\]

Now we will compute \ \m{\Ext^1_{\ko_\kc}(\ko_{I^c\backslash\{n\}},
\ko_{\{n\}})}. According to \cite{dr4}, 2.3, we have an exact sequence
\[0\lra\Ext^1_{\ko_{I^c}}(\ko_{I^c\backslash\{n\}},\ko_{\{n\}})\lra
\Ext^1_{\ko_\kc}(\ko_{I^c\backslash\{n\}},\ko_{\{n\}})\lra\Hom(
\Tor^1_{\ko_\kc}(\ko_{I^c\backslash\{n\}},\ko_{I^c}),\ko_{\{n\}}) .\]
Since \m{\Tor^1_{\ko_\kc}(\ko_{I^c\backslash\{n\}},\ko_{I^c})} is concentrated
on \m{\kc_{I^c\backslash\{n\}}}, we have
\[\Hom(\Tor^1_{\ko_\kc}(\ko_{I^c\backslash\{n\}},\ko_{I^c}),\ko_{\{n\}})
\ = \ \nsp .\]
So we have
\[\Ext^1_{\ko_\kc}(\ko_{I^c\backslash\{n\}},\ko_{\{n\}}) \ = \
\Ext^1_{\ko_{I^c}}(\ko_{I^c\backslash\{n\}},\ko_{\{n\}}) .\]
Let \m{\kj} denote the ideal sheaf of \m{\kc_{\{n\}}} in \m{\kc_{I^c}}. The
ideal sheaf of \m{\kc_{I^c\backslash\{n\}}} is generated by
\m{{\bf w}=(0,\ldots,0,\pi_n^m)}, with \ \m{m=\sigg_{i\in
I^c\backslash\{n\}}p_{in}}. So we have an exact sequence of sheaves on
\m{\kc_{I^c}}
\xmat{ 0\ar[r] & \kj\ar[r] & \ko_{I^c}\ar[r]^-\alpha & \ko_{I^c}\ar[r] &
\ko_{I^c\backslash\{n\}}\ar[r] & 0 \ ,}
where $\alpha$ is the multiplication by \m{\bf w}. By the induction hypothesis
there exists a surjective morphism \m{\ko_{I^c}\to\kj}, so we get a locally
free resolution of \m{\ko_{I^c\backslash\{n\}}}
\xmat{\ko_{I^c}\ar[r] & \ko_{I^c}\ar[r]^-\alpha & \ko_{I^c}\ar[r] &
\ko_{I^c\backslash\{n\}}\ar[r] & 0 \ ,}
that can be used to compute
\m{\EExt^1_{\ko_{I^c}}(\ko_{I^c\backslash\{n\}},\ko_{\{n\}})}. It follows
easily that
\[\EExt^1_{\ko_{I^c}}(\ko_{I^c\backslash\{n\}},\ko_{\{n\}}) \ \simeq \
\ko_{\{n\}}/(\pi_n^m) \ .\]
We have \ \m{\HHom(\ko_{I^c\backslash\{n\}},\ko_{\{n\}})=0}, hence
\begin{eqnarray*}
\Ext^1_{\ko_{I^c}}(\ko_{I^c\backslash\{n\}},\ko_{\{n\}}) & \simeq & H^0(
\EExt^1_{\ko_{I^c}}(\ko_{I^c\backslash\{n\}},\ko_{\{n\}}))\\
& \simeq & H^0(\ko_{\{n\}}/(\pi_n^m))\\
& \simeq & H^0(\ko_S/(\pi_n^m))\\
& \simeq & \C[\pi_n]/(\pi_n^m) .
\end{eqnarray*}
We will now describe the sheaves $\ke$ such that there exists an exact
sequence
\begin{equation}\label{ecl20}
0\lra\ko_{\{n\}}\lra\ke\lra\ko_{I^c\backslash\{n\}}\lra 0 .
\end{equation}
Let \ \m{\nu\in \C[\pi_n]/(\pi_n^m)} \ be associated to this exact sequence,
and \ \m{\ov{\nu}\in H^0(\ko_S)} \ over $\nu$. Let
\[{\xymatrix@R=4pt{
\tau : \ko_{\{n\}}\ar[r] & \ko_{\{n\}}\oplus\ko_{I^c} \\
\ \ \ \ u\fmaps[r] & (\ov{\nu} u,{\bf w}u) \\
}}\]
Then according to the preceding resolution of \m{\ko_{I^c\backslash\{n\}}} and
the construction of extensions (cf. \cite{dr3}, 4.2), we have \
\m{\ke\simeq\coker(\tau)}. It is easy to see that if \m{\nu=-1} then \
\m{\ke\simeq\ko_{\ki^c}}. If $\nu$ is invertible, then we have also \
\m{\ke\simeq\ko_{\ki^c}}, because the corresponding extension can be obtained
from the one corresponding to \m{\nu=-1} by multiplying the left morphism of
the exact sequence by $\nu$.

A similar construction can be done for extensions of
\m{\ko_{I^c,x}}-modules (for every \m{x\in C})
\[0\lra\ko_{\{n\},x}\lra V \lra\ko_{I^c\backslash\{n\},x}\lra 0 .\]
These extensions are classified by \m{\ko_{\{n\},x}/(\pi_n^m)}, and
\m{\ko_{I^c,x}} corresponds to $-1$.

Conversely we consider extensions
\xmat{0\ar[r] & \ko_{\{n\},x} \ar[r]^-\lambda & \ko_{I^c,x}\ar[r]^-\mu &
\ko_{I^c\backslash\{n\},x}\ar[r] & 0 \ .}
Using the facts that \m{\Hom(\ko_{\{n\},x},\ko_{I^c,x})} is generated by
the multiplication by $\bf w$ and \Nligne
\m{\Hom(\ko_{I^c,x},\ko_{I^c\backslash\{n\},x})} by the restriction morphism, 
it is easy to see that $\lambda$, $\mu$ are unique up to multiplication by an
invertible element of \m{\ko_{I^c,x}}. Hence the elements of
\m{\Ext^1_{\ko_{I^c,x}}(\ko_{I^c\backslash\{n\},x},\ko_{\{n\},x})} 
corresponding to the preceding extensions are exactly the invertible elements 
of \m{\ko_{\{n\},x}/(\pi_n^m)}.

It follows that the extensions \eqref{ecl20} where $\ke$ is locally free
correspond to invertible elements of \m{\C[\pi_n]/(\pi_n^m)}, and we have seen
that in this case we have \m{\ke\simeq\ko_{I^c}}. Hence we have \
\m{\ki_{\kc_I}\simeq\ko_{I^c}} \ and theorem \ref{ecl19} is proved.
\end{proof}

\sepprop

\begin{subsub}\label{ecl21}{\bf Corollary: } The ideal sheaf of \m{\kc_I} is
globally generated by an element \m{{\bf u}_I} such that for every integer $i$
such that \m{1\leq i\leq n} and \m{i\not\in I}, the $i$-th coordinate of
\m{{\bf u}_I} belongs to \m{H^0(\ko_S)}.
\end{subsub}

\end{sub}

\sepsub

\Ssect{Properties of the fragmented deformations}{const_def}

We use the notations of \ref{spectre}.

Let $i$ be an integer such that \m{1\leq i\leq n} and
\m{J_i=\{1,\ldots,n\}\backslash\{i\}}.
We denote by \m{\kb} the image of \m{\ko_{\kc}} in \
\m{\prod_{1\leq j\leq n}\ko_{\kc_j}/(\pi_j^{q_j})}; it is a sheaf of
$\C$-algebras on $C$.
Let \m{\kb_i} be the image of
\m{\ko_{\kc_{J_i}}} in \ \m{\prod_{1\leq j\leq
n,j\not=i}\ko_{\kc_j}/(\pi_j^{q_j})}; it is also a sheaf of $\C$-algebras on
$C$. For every point $x$ of $C$ and every
\m{\alpha=(\alpha_m)_{1\leq m\leq n}} in \ \m{\prod_{1\leq j\leq 
n}\ko_{\kc_j,x}}, we denote by \m{b_{i}(\alpha)} its image in \
\m{\prod_{1\leq j\leq n,j\not=i}\ko_{\kc_j,x}} (obtained by
forgetting the $i$-th coordinate of $\alpha$).

If \m{p,k} are positive integers, with \m{k\leq n}, \m{x\in C} and 
\m{\alpha\in\ko_{\kc_k,x}}, let \m{[\alpha]_p} denote the image of $\alpha$ in 
\m{\ko_{\kc_k,x}/\pi_k^p}.

\sepprop

\begin{subsub}\label{ecl9}{\bf Proposition: } There exists a morphism of
sheaves of algebras on $C$
\[\Phi_i : \kb_i\lra\ko_{\kc_i}/(\pi_i^{q_i})\]
such that for every point $x$ of $C$ and all \m{(\alpha_m)_{1\leq m\leq
n,m\not=i}\in\ko_{\kc_{J_i},x}}, \m{\alpha_i\in\ko_{\kc_i,x}}, we have \Nligne
\m{\alpha=(\alpha_m)_{1\leq m\leq n}\in\ko_{\kc,x}} if and only if \
\m{\Phi_{i,x}(b_{i}(\alpha))=[\alpha_i]_{q_i}}.
\end{subsub}

\begin{proof}
Let \ \m{(\alpha_m)_{1\leq m\leq n,m\not=i}\in\ko_{\kc_{J_i}x}}. Since
\m{\kc_{J_i}\subset\kc}, there exists \m{\alpha_i\in\ko_{\kc_i,x}} such that
\Nligne
\m{(\alpha_m)_{1\leq m\leq n}\in\ko_{\kc,x}}. If \m{\alpha'_i\in\ko_{\kc_i}x}
has the same property, we have
\Nligne\m{(0,\ldots,0,\alpha_i-\alpha'_i,0,\ldots,0)\in\ki_{J_ix}}. So
according to proposition \ref{ecl8b}, we have \
\m{[\alpha_i]_{q_i}=[\alpha'_i]_{q_i}}. Hence we have a well defined morphism
of algebras \ \m{\theta_x:\ko_{\kc_{J_i},x}\to\ko_{\kc_{J_i}}/(\pi_i^{q_i})} \
sending \m{(\alpha_m)_{1\leq m\leq n,m\not=i}} \ to \m{[\alpha_i]_{q_i}}. If
\m{j\in J_i}, we have, according to proposition \ref{ecl8b},\Nligne
\m{\theta_x(0,\ldots,0,\pi_j^{q_j},0,\ldots,0)=0}. Hence $\theta_x$ induces a
morphism of algebras \m{\kb_{i,x}\to\ko_{\kc_i,x}/(\pi_i^{q_i})}.
\end{proof}

\sepprop

The morphism \m{\Phi_i} has the following properties: for every point $x$ of
$C$
\begin{enumerate}
\item[(i)] For every \ $\alpha=(\alpha_m)_{1\leq m\leq n,m\not=i}\in\kb_{i,x}$,
we have \ $\Phi_{i,x}(\alpha)_{\mid C}=\alpha_{m\mid C}$ \ for $1\leq m\leq 
n$, $m\not=i$.
\item[(ii)] We have \ $\Phi_{i,x}((\pi_m)_{1\leq m\leq n,m\not=i})=\pi_i$.
\item[(iii)] Let \ $j,k\in\{1,\cdots,n\}$ \ be such that $i,j,k$ are distinct.
Let $\bf v$ be the image of ${\bf u}_{jk}$ in $\kb_i$.
Then there exists $\lambda\in\ko_{\kc_i,x}^*$ such that \
$\Phi_{i,x}({\bf v})=\lambda\pi_i^{p_{ij}}$.
\item[(iv)] Let $j$ be an integer such that $1\leq j\leq n$ and $j\not=i$. Let
$\bf v$ be the image of ${\bf u}_{ij}$ in $\kb_{i,x}$. Then we have \
$\ker(\Phi_{i,x})=({\bf v})$.
\end{enumerate}

\sepprop

\begin{subsub}\label{recip} Converse - \rm Let $\kc'$ be a gluing of
\m{\kc_1,\ldots,\kc_{i-1},\kc_{i+1},\ldots,\kc_n} along $C$, which is a 
fragmented 
deformation of a primitive multiple curve of multiplicity \m{n-1}.  Let
\m{(p_{jk})_{1\leq j,k\leq n,j,k\not=i}} be the spectrum of \m{\kc'}. 
Let \m{p_{ij}}, \m{1\leq j\leq n,j\not=i} be positive integers, and 
\m{p_{ii}=0}. For \m{1\leq j\leq n}, let \m{q_j=\sigg_{1\leq k\leq n}p_{kj}}.

Let \m{\kb_i} be the image of \m{\ko_{\kc'}} in \ \m{\prod_{1\leq j\leq
n,j\not=i}\ko_{\kc_j}/(\pi_j^{q_j})} \ and
\[\Phi_i : \kb_i\lra\ko_{\kc_i}/(\pi_i^{q_i})\]
a morphism of sheaves of algebras on $C$ satisfying properties
(i), (ii), (iii) above. Let $\ka$ be the subsheaf of algebras of $\BS{\ka}$
defined by: \m{\ka=\BS{\ka}} on \m{\BS{\ka}_{top}\backslash C}, and for every
point $x$ of $C$, and every \ \m{\alpha=(\alpha_m)_{1\leq m\leq
n}\in\prod_{m=1}^n\ko_{\kc_{m,x}}}, \m{\alpha\in\ka_x} if and only if \ 
\m{b_i(\alpha)\in\kb_{i,x}} \ and
\m{\Phi_{i,x}(b_i(\alpha))=\lbrack\alpha_i\rbrack_{q_i}}.

It is easy to see that $\ka$ is the structural sheaf of a gluing of
\m{C_1,\ldots,C_n} along $C$, which is a fragmented deformation of a primitive 
multiple curve of multiplicity $n$, and that \ 
\m{\kc'=\ka_{\{1,\ldots,i-1,i+1,\ldots,n\}}}.
\end{subsub}
\sepprop

We give now some applications of the preceding construction.

\sepprop

\begin{subsub}\label{coro1}{\bf Corollary: } Let $N$ be an integer such that
\m{N\geq\max_{1\leq i\leq n}(q_i)}. Let \m{x\in C},\Nligne
\m{\beta\in\ko_{\kc_1,x}\times\cdots\ko_{\kc_n,x}} \ and \m{u\in\ko_{\kc,x}} 
such
that \ \m{u_{\mid C}\not=0}. Suppose that \ \m{\lbrack\beta
u\rbrack_N\in\ko_{\kc,x}/(\pi^N)}. Then we have \
\m{\lbrack\beta\rbrack_N\in\ko_{\kc,x}/(\pi^N)}.
\end{subsub}

\begin{proof} By induction on $n$. It is obvious if $n=1$. Suppose that the
lemma is true for \m{n-1}. Let \ \m{I=\{1,\ldots,n-1\}}. So we have \
\m{[\beta_{\mid\kc_1\times\cdots\kc_{n-1}}]_N 
\in\ko_{\kc_I,x}/(\pi_1,\ldots,\pi_{n-1})^N} by the induction hypothesis.
Let $\gamma$ (resp. $v$) be the image of \m{\beta} (resp. $u$) in \m{\kb_n}. To
show that \m{\lbrack\beta\rbrack_N\in\ko_{\kc,x}/(\pi^N)} it is enough to
verify that
\[\Phi_n(\gamma) \ = \ \lbrack\beta_n\rbrack_{q_n} .\]
We have \ \m{\Phi_n(\gamma v) \ = \ \lbrack\beta_nu_n\rbrack_{q_n}} because
\m{\lbrack\beta u\rbrack_N\in\ko_{\kc,x}/(\pi^N)}, and \m{\Phi_n(v)=\lbrack u_n
\rbrack_{q_n}} because \ \m{u\in\ko_{\kc,x}}. So we have
\[\Phi_n(\gamma)\lbrack u_n\rbrack_{q_n} \ = \ \Phi_n(\gamma)\Phi_n(v) \ = \
\Phi_n(\gamma v) \ = \ \lbrack\beta_nu_n\rbrack_{q_n} \ = \
\lbrack\beta_n\rbrack_{q_n}\lbrack u_n\rbrack_{q_n} .\]
Since \m{u_{\mid C}\not=0}, \m{\lbrack u_n\rbrack_{q_n}} is not a zero divisor
in \m{\ko_{\kc_n,x}/(\pi_n^{q_n})}, so we have \ \m{\Phi_n(\gamma)=
\lbrack\beta_n\rbrack_{q_n}}.
\end{proof}

\sepprop

\begin{subsub}\label{coro2}{\bf Corollary: } Let \ \m{{\bf q}=\max_{1\leq
i\leq n}(q_i)} \ and $p$ the number of integers $i$ such that \m{1\leq i\leq
n} and \ \m{q_i={\bf q}}. Then we have \ \m{p\geq 2}.
\end{subsub}

\begin{proof} Suppose that \ \m{q_i={\bf q}}. Then we have \
\m{\pi_i^{q_i-1}\not=0} \ in \m{\ko_{\kc_i}/(\pi_i^{q_i})}. Since \Nligne
\m{\pi_i=\Phi_i((\pi_m)_{1\leq m\leq n,m\not=i})}, we have \
\m{(\pi_m^{q_i-1})_{1\leq m\leq n,m\not=i}\not=0} \ in \m{\kb_i}. So we
cannot have  \ \m{q_m<q_i} \ for all the \m{m\not=i}.
\end{proof}

Let $i$ be an integer such that \m{1\leq i\leq n},
\[\kh=\prod_{1\leq j\leq n}(\pi_j^{q_j-1})/(\pi_j^{q_j})\simeq\ko_C^n \quad
\quad\quad \text{(resp.} \quad\kh_i=\prod_{1\leq j\leq n,j\not=i}
(\pi_j^{q_j-1})/(\pi_j^{q_j})\simeq\ko_C^{n-1} \ ) .\]
It is an ideal sheaf of \ \m{\prod_{1\leq j\leq n}\ko_{\kc_j}/(\pi_j^{q_j})}
(resp. \m{\prod_{1\leq j\leq n,j\not=i}\ko_{\kc_j}/(\pi_j^{q_j})} ). Let \
\m{\kj=\kh\cap\kb} (resp. \m{\kj_i=\kh_i\cap\kb_i}), which is an ideal sheaf
of $\kb$ (resp. \m{\kb_i}).

\sepprop

\begin{subsub}\label{ecl10} {\bf Proposition: } There exists a unique \
\m{\lambda(\kc)=(\lambda_1,\ldots,\lambda_n)\in\P_n(\C)} such that for
every\Nligne \m{{\bf u}=(u_j)_{1\leq j\leq n}\in\kh}, we have \m{{\bf
u}\in\kj} if and only if \ \m{\lambda_1u_1+\cdots\lambda_nu_n=0}. The
\m{\lambda_i} are all non zero.
\end{subsub}

\begin{proof} We have \ \m{(\pi_m)_{1\leq m\leq n,m\not=i}.\kj_i=0}. Hence \
\m{\pi_i\Phi_i(\kj_i)=0} \ and \Nligne \m{\Phi_i(\kj_i)\subset(\pi_i^{q_i-1})/
(\pi_j^{q_i})}. The restriction of \m{\Phi_i}, \m{\kj_i\to(\pi_i^{q_i-1})/
(\pi_j^{q_i})} \ is a morphism \Nligne \m{(n-1)\ko_C\to\ko_C} \ of vector
bundles on $C$. The existence of \m{(\lambda_1,\ldots,\lambda_n)} follows from
that.

If \m{\lambda_i=0}, we have \ \m{(0,\ldots,0,\pi_i^{q_i-1},0,\ldots,0)\in
\ko_\kc(\kc)}. This is impossible because according to proposition \ref{ecl8b},
\m{(0,\ldots,0,\pi_i^{q_i},0,\ldots,0)} generates the ideal sheaf of
\m{\kc_{J_i}} in $\kc$.
\end{proof}

\sepprop

For all distinct integers $i$, $j$ such that \m{1\leq i,j\leq n}, let \
\m{I_{ij}=\{1,\ldots,n\}\backslash\{i,j\}}. Then according to proposition
\ref{ecl8b}, \m{{\bf u}_{I_{ij}i}} generates the ideal sheaf of
\m{\kc_{I_{ij}}}. We have \ \m{{\bf u}_{I_{ij}i}=(b_k)_{1\leq k\leq n}}, with \
\m{b_k=0} \ if \m{k\not=i,j}, \m{b_i=\pi_i^{q_i-p_{ij}}} \ and
\[b_j \ = \ \big(\prod_{1\leq m\leq n,m\not=i,j}\alpha^{(j)}_{mi} \ \big).
\pi_j^{q_j-p_{ij}} .\]
So we have \ \m{\pi^{p_{ij}-1}{\bf u}_{I_{ij}i}\in\kj_i}, which gives the
equation
\begin{equation}\label{equ2}
\frac{\lambda_i}{\lambda_j} \ = \ - \prod_{1\leq m\leq
n,m\not=i,j}{\bf a}^{(j)}_{mi} .
\end{equation}

\sepprop

\begin{subsub}\label{ecl11}{\bf Proposition: } For all distinct integers
$i,j,k$ such that \m{1\leq i,j,k\leq n}, we have
\[{\bf a}_{ki}^{(j)} \ = \ - {\bf a}_{ik}^{(j)}{\bf a}_{ji}^{(k)} \ . \]
\end{subsub}

\begin{proof} We need only to treat the case \m{n=3}, and we get the preceding
formula by writing that \ \m{\frac{\lambda_1}{\lambda_3}=
\frac{\lambda_1}{\lambda_2}.\frac{\lambda_2}{\lambda_3}}, and by using 
\eqref{equ2}.
\end{proof}

\sepprop

\begin{subsub}\label{ecl14}{\bf Proposition: } Let \ \m{(\alpha_1\pi_1^{m_1},
\ldots,\alpha_n\pi_n^{m_n})\in\ko_{\kc,x}}, with \m{\alpha_1,\ldots,\alpha_n}
invertible. Let \ \m{M=m_1+\cdots+m_n}. then
\[\big(\frac{1}{\alpha_1}\pi_1^{M-m_1},\ldots,\frac{1}{\alpha_n}\pi_n^{M-m_n}
\big) \ \in \ \ko_{\kc,x} \ .\]
\end{subsub}

\begin{proof} By induction on $n$. It is obvious for \m{n=1}. Suppose that it
is true for \m{n-1\geq 1}. Let \ \m{I=\{1,\ldots,n-1\}}. Then \
\m{(\alpha_1\pi_1^{m_1},\ldots,\alpha_{n-1}\pi_{n-1}^{m_{n-1}})\in
\ko_{\kc_{I,x}}}. Hence, by the induction hypothesis, we have
\[\big(\frac{1}{\alpha_1}\pi_1^{M-m_1-m_n},\ldots,
\frac{1}{\alpha_{n-1}}\pi_{n-1}^{M-m_{n-1}-m_n}
\big) \ \in \ \ko_{\kc_{I,x}} \ .\]
So there exists \m{\gamma\in\ko_{\kc_{n,x}}} such that
\[u \ = \ \big(\frac{1}{\alpha_1}\pi_1^{M-m_1-m_n},\ldots,
\frac{1}{\alpha_{n-1}}\pi_{n-1}^{M-m_{n-1}-m_n},\gamma\big)
\ \in \ \ko_{\kc,x} \ .\]
Multiplying by  \m{(\alpha_1\pi_1^{m_1},\ldots,\alpha_n\pi_n^{m_n})} we see
that \ \m{(\pi_1^{M-m_n},\ldots,\pi_{n-1}^{M-m_n},\gamma
\alpha_n\pi_n^{m_n})\in\ko_{\kc,x}}. Subtracting \m{\pi^{M-m_n}}, we find
that \ \m{(0,\ldots,0,\gamma\alpha_n\pi_n^{m_n}-\pi_n^{M-m_n})\in\ko_{\kc,x}}.
There exists \m{\alpha\in\ko_{\kc,x}} such that the $n$-th coordinate of
$\alpha$ is \m{\alpha_n}, and $\alpha$ is invertible. It follows that \Nligne
\m{v=(0,\ldots,0,\gamma\pi_n^{m_n}-\frac{1}{\alpha_n}\pi_n^{M-m_n})\in
\ko_{\kc,x}}. Now we have
\[\pi^{m_n}u-v \ = \
\big(\frac{1}{\alpha_1}\pi_1^{M-m_1},\ldots,\frac{1}{\alpha_n}\pi_n^{M-m_n}
\big) \ \in \ \ko_{\kc,x} \ . \]
\end{proof}

\sepprop

\begin{subsub}\label{ecl22}{\bf Corollary: } Let \m{V\subset U} be open
subsets of $\kc$, and suppose that \ \m{U\cap C\not=\emptyset}. Let
\m{\alpha\in\ko_{\kc}(V)} and \m{\beta\in\ko_{\BS{\ka}}(U)} such that \
\m{\beta_{\mid V}=\alpha}. Then \ \m{\beta\in\ko_{\kc}(U)}.
\end{subsub}

(Recall that $\BS{\ka}$ is the {\em initial gluing} of \m{\kc_1,\ldots,\kc_n}
(cf. \ref{ecl2c})).

\begin{proof} This can be proved easily by induction on $n$, using proposition
\ref{ecl9}.
\end{proof}

\end{sub}

\sepsub

\Ssect{Construction of fragmented deformations}{const_def2}

Consider a fragmented deformation
\[\pi=\pi^{[n-1]}=(\pi_1,\ldots,\pi_{n-1}) : \kc^{[n-1]}\lra S\]
of \m{C_{n-1}}, with \m{n-1} irreducible components
\m{\kc_1,\ldots,\kc_{n-1}}. Let \m{(p_{ij}^{[n-1]})_{1\leq i,j<n}} be its
spectrum. For \m{1\leq i<n}, let \ \m{q_i^{[n-1]}=\sigg_{1\leq
j<n}p_{ij}^{[n-1]}}. We denote by \m{\ki_C^{[n-1]}} the ideal sheaf of $C$ in
\m{\kc^{[n-1]}}. Let \
\m{\lambda(\kc^{[n-1]})=(\lambda_1,\ldots,\lambda_{n-1})}.

Let \m{p_{1n},\ldots,p_{n-1,n}} be positive integers,
\m{q_i=q_i^{[n-1]}+p_{in}} \ for \m{1\leq i<n}, and \Nligne
\m{q_n=p_{1n}+\cdots+p_{n-1,n}}.
Let \m{{\bf u}\in H^0(\ki_{C}^{[n-1]})} whose image
generates \m{\ki_{C}^{[n-1]}/((\ki_{C}^{[n-1]})^2+(\pi))}, of the form
\[{\bf u} \ = \ (\beta_1\pi_1^{p_{1n}},\ldots,\beta_{n-1}
\pi_{n-1}^{p_{n-1,n}}) ,\]
with \m{\beta_i\in H^0(\ko_S)} invertible for \m{1\leq i<n}.

Let \m{\kb^{[n-1]}} be the image of \m{\ko_{\kc^{[n-1]}}} in \ \m{\ko_{\kc_1}/
(\pi_1^{q_1})\times\cdots\times\ko_{\kc_{n-1}}/ (\pi_{n-1}^{q_{n-1}})} . We
will also denote by $\bf u$ the image of $\bf u$ in \m{\kb^{[n-1]}}. Let \
\m{\kq=\kb^{[n-1]}/({\bf u})} , \m{\rho:\kb^{[n-1]}\to\kq} \ the projection and
\ \m{\pi_n=\rho(\pi)}.

\sepprop

\begin{subsub}\label{ecl15}{\bf Proposition: } We have \ \m{\pi_n^{q_n}=0} .
\end{subsub}

\begin{proof} According to proposition \ref{ecl14} we have for every \m{x\in C}
\[v \ = \ \big(\frac{1}{\beta_1}\pi_1^{q_n-p_{1n}},\ldots,\frac{1}{\beta_{n-1}}
\pi_n^{q_n-p_{n-1,n}}\big) \ \in \ \ko_{\kc^{[n-1]} x} \ .\]
Hence \ \m{\pi^{q_n}=v{\bf u}\in({\bf u})} \ in \m{\ko_{\kc^{[n-1]} x}}, and
\ \m{\pi_n^{q_n}=0} .
\end{proof}

\sepprop

Now we suppose that we have \m{p_{ij}=p_{kl}}, for \m{1\leq i,j,k,l\leq n}, 
\m{i\not=j}, \m{k\not=l}. This is the case for example if all the components 
\m{\kc_i} are transverse in $\kc$ (i.e. \m{p_{ij}=1}, \m{1\leq i,j\leq 
n}, \m{i\not=j}).

\sepprop

\begin{subsub}\label{ecl17}{\bf Proposition: 1 --} We have \
\m{\pi_n^{q_n-1}=0} \ if and only if
\[\frac{\lambda_1}{\beta_{1\mid C}}+\cdots+
\frac{\lambda_{n-1}}{\beta_{n-1\mid C}} \ = \ 0 \ .\]

We suppose now that \ \m{\frac{\lambda_1}{\beta_{1\mid C}}+\cdots+
\frac{\lambda_{n-1}}{\beta_{n-1\mid C}}\not=0}. Let \m{x\in C}. Then

{\bf 2 --} For every \ \m{\epsilon\in\kb^{[n-1]}_x} \ such that \
\m{\epsilon_{\mid C}\not=0}, we have \ \m{\pi^{q_n-1}\epsilon\not\in({\bf u})}.

{\bf 3 --} For every \ \m{\eta\in\kb^{[n-1]}_x/({\bf u)}}, and every integer
$k$ such that \m{1\leq k<q_n}, we have \ \m{\pi_n^k\eta=0} \ if and only if
$\eta$ is a multiple of \m{\pi_n^{q_n-k}}.

{\bf 4 --} \m{\kb^{[n-1]}_x/({\bf u)}} is a flat
\m{\C[\pi_n]/(\pi_n^{q_n})}-module.
\end{subsub}

\begin{proof}
We have \m{\pi_n^{q_n-1}=0} if and only if \
\m{(\pi_1^{q_n-1},\ldots,\pi_{n-1}^{q_n-1})\in({\bf u})} \ in
\m{\kb_x^{[n-1]}}.
We have, in \m{\ko_{\kc_{1x}}\times\cdots\times\ko_{\kc_{n-1,x}}},
\[(\pi_1^{q_n-1},\ldots,\pi_{n-1}^{q_n-1}) \ = \
(\beta_1\pi_1^{p_{1n}},\ldots,\beta_{n-1,n}\pi_{n-1}^{p_{n-1,n}}).
(\frac{1}{\beta_1}\pi_1^{q_1^{[n-1]}-1},\ldots,
\frac{1}{\beta_p}\pi_1^{q_{n-1}^{[n-1]}-1}) ,\]
and \m{\pi_n^{q_n-1}=0} if and only if there exist
\m{\eta\in\ko_{\kc^{[n-1]},x}}, \m{a_i\in\ko_{\kc_i,x}}, \m{1\leq i<n}, such 
that
\[(\pi_1^{q_n-1},\ldots,\pi_{n-1}^{q_n-1}) \ = \ \eta{\bf u}+
(a_1\pi_1^{q_1},\ldots,a_{n-1}\pi_{n-1}^{q_{n-1}}) . \]
This equality is equivalent to
\[(\frac{1}{\beta_1}\pi_1^{q_1^{[n-1]}-1},\ldots,
\frac{1}{\beta_{n-1}}\pi_1^{q_{n-1}^{[n-1]}-1}) -\eta \ = \
(\frac{a_1}{\beta_1}\pi_1^{q_1^{[n-1]}},\ldots,
\frac{a_{n-1}}{\beta_{n-1}}\pi_{n-1}^{q_{n-1}^{[n-1]}}) . \]
Since for \m{1\leq i<n}, we have \ \m{(0,\ldots,0,\pi_i^{q_i^{[n-1]}},0,\ldots
0)\in\ko_{\kc^{[n-1]},x}}, we have \m{\pi_n^{q_n-1}=0} if and only if
\[(\frac{1}{\beta_1}\pi_1^{q_1^{[n-1]}-1},\ldots,
\frac{1}{\beta_{n-1}}\pi_1^{q_{n-1}^{[n-1]}-1}) \ \in \ \ko_{\kc^{[n-1]},x}
. \]
So the result of 1- follows from the definition of \m{\lambda(\kc^{[n-1]})}
(cf. prop. \ref{ecl10}), 2- is an easy consequence.

Now we prove 3-, by induction on $k$. Suppose that it is true for \m{k=1}, and
that \m{\pi^k_n\eta=0}, with \m{2\leq k<q_n}. We have
\m{\pi^{k-1}_n.\pi_n\eta=0}, so according to the induction hypothesis,
\m{\pi_n\eta} is a multiple of \m{\pi_n^{q_n-k+1}}: \m{\pi_n\eta=
\pi_n^{q_n-k+1}\lambda}. So \ \m{\pi_n(\eta-\pi_n^{q_n-k}\lambda)=0}.
Since -3 is true for \m{k=1}, we can write \
\m{\eta-\pi_n^{q_n-k}\lambda=\pi_n^{q_n-1}\epsilon}, i.e. \ \m{\eta
=\pi_n^{q_n-k}(\lambda+\pi_n^{k-1}\epsilon)}, and 3- is true for $k$.

Il remains to prove 3- for \m{k=1}. Suppose that \ \m{\pi_n\eta=0} (with
\m{\eta\not=0}). We can write $\eta$ as \ \m{\eta=\pi_n^m\theta}, where
$\theta$ is not a multiple of $\pi_n$, and \m{0\leq m<q_n}. Let \
\m{\ov{\theta}\in\kb^{[n-1]}_x} \ be over $\theta$. Since \ 
\m{\ki_C^{[n-1]}=({\bf
u})+(\pi)} \ according to proposition \ref{ecl12}, the condition ``$\theta$ is
not a multiple of $\pi_n$'' is equivalent to \ 
\m{\ov{\theta}\not\in\ki_{C,x}^{[n-1]}}.
We have \ \m{\pi^{m+1}\ov{\theta}\in({\bf u})}, so according to 2-, we have \
\m{m+1\geq q_n}, which proves 3- for \m{k=1}. The last assertion is an easy
consequence of 3-.
\end{proof}

\sepprop

\begin{subsub}\label{exe1}{\bf Example: } \rm Let $N$ be an integer,
\m{s\in H^0(\ko_S)} invertible, and \m{k,l} integers such that
\m{1\leq k,l<n}, \m{k\not=l}. Suppose that for every integer $i$ such that
\m{1\leq i<n} and \m{i\not=k} we have \ \m{N>p_{ik}^{[n-1]}} \ and \ \m{N\geq
q_i^{[n-1]}-q_k^{[n-1]}+p_{ik}^{[n-1]}}. We take \ \m{{\bf u}={\bf
u}_{kl}-s\pi^N}. We have then \ \m{\beta_i=\alpha^{(i)}_{kl}} \ if \m{i\not=k},
and \ \m{\beta_k=-s}. The condition \ \m{\frac{\lambda_1}{\beta_{1\mid
C}}+\cdots+\frac{\lambda_{n-1}}{\beta_{n-1\mid C}}\not=0} \ is fulfilled
if and only if
\[ \sigg_{1\leq i<n,i\not=k}
\frac{\lambda_i}{{\bf a}_{kl}^{(i)}}-\frac{\lambda_k}{s_{\mid C}} \ \not= \ 0 
. \]
\end{subsub}

\sepprop

\begin{subsub}\label{cst_fd} Construction of fragmented deformations -- \rm
Suppose that \ \m{\frac{\lambda_1}{\beta_{1\mid C}}+\cdots+
\frac{\lambda_{n-1}}{\beta_{n-1\mid C}}\not=0}. From proposition \ref{ecl17},
4-, it is easy to prove that
\begin{enumerate}
\item[--] There exists a flat morphism of algebraic varieties \
\m{\tau:Y\to\spec(\C[\pi_n]/(\pi_n^{q_n}))} \ with a canonical isomorphism of
sheaves of \m{\C[\pi_n]/(\pi_n^{q_n})}-algebras \ \m{\ko_Y\simeq\kq}, such
that \ \m{\tau^{-1}(*)=C} (where $*$ is the closed point of
\m{\spec(\C[\pi_n]/(\pi_n^{q_n}))}).
\item[--] There exist a family of smooth curves \m{\kc_n} and a flat
morphism \ \m{\pi_n:\kc_n\to S} \ extending $\tau$ (recall that $S$ is a
germ). Hence $Y$ is the inverse image of the subscheme of \m{\kc_n}
corresponding to the ideal sheaf \m{(\pi_n^{q_n})}. The existence of \m{\kc_n}
can be proved using Hilbert schemes of curves in projective spaces. Of course
\m{\kc_n} need not be unique.
\end{enumerate}
We obtain a gluing $\kc$ of \m{\kc_1,\ldots,\kc_n} by defining the sheaves of
algebras \m{\ko_\kc} (on the Zariski topological space corresponding to the
initial gluing $\BS{\ka}$) as in \ref{recip}, using for \m{\Phi_n} the
quotient morphism \m{\kb^{[n-1]}\to\kq}. It is easy to see that
\m{\pi^{-1}(P)} is a primitive multiple curve \m{C_n} of multiplicity $n$
extending \m{C_{n-1}}, hence $\kc$ is a fragmented deformation of $C_n$.
\end{subsub}

\sepprop

\begin{subsub}\label{rem_fd}{\bf Remark: }\rm {\bf 1 -- } The multiple curve
\m{C_n} depends on the choice of the family \m{\kc_n} extending the family $Y$
parametrized by \m{\spec(\C[\pi_n]/(\pi_n^{q_n}))}.

{\bf 2 -- } The multiple curve \m{C_{n-1}} is completely defined by
\m{\kb^{[n-1]}}, because \ \m{(\pi_1^{q_1})\times\cdots(\pi_{n-1}^{q_{n-1}})
\subset(\pi)}. But it is not enough to know  \m{\kb^{[n-1]}} and \m{\bf u} to
define \m{C_n}. In fact we need \m{\ko_{\kc_i}/(\pi_i^{q_i+1})}, \m{1\leq
i\leq n}.
\end{subsub}

\end{sub}

\sepsub

\Ssect{Basic elements}{pol_elem}

We use the notations of \ref{spectre} and \ref{const_def}.

Let \m{{\bf m}=(m_1,\ldots,m_n)} be an $n$-tuple of positive integers, and
\[\BS{\Pi}^{\bf m} \ = \ (\pi_1^{m_1})\times\cdots\times(\pi_n^{m_n}) .\]

\sepprop

\begin{subsub}\label{def0}{\bf Definition: } Let \m{x\in C}. An element $u$
of \m{\ko_{\kc,x}} is called {\em basic} {\em at order} {\bf m} if there
exist polynomials \m{P_1,\ldots,P_n\in\C[X]} such that
\[u \ \equiv \ (P_1(\pi_1),\ldots,P_n(\pi_n)) \quad\quad (\text{\em mod.}
\quad \BS{\Pi}^{\bf m}) \ .\]
If \m{u=(P_1(\pi_1),\ldots,P_n(\pi_n))}, we say that $u$ is {\em basic}.
\end{subsub}

\sepprop

Let \ \m{{\bf q}=(q_1,\ldots,q_n)}. Then according to corollary \ref{ecl22},
if $u$ is basic at order $\bf q$, then for every \m{y\in C}, we
have \m{(P_1(\pi_1),\ldots,P_n(\pi_n))\in\ko_{\kc,y}}. So
\m{(P_1(\pi_1),\ldots,P_n(\pi_n))} is defined on a neighborhood of $C$.

\sepprop

\begin{subsub}\label{ecl23}{\bf Lemma: } Let \m{u,v,w\in\ko_{\kc,x}}
such that \ \m{w=uv} \ and \ \m{w\not=0}. Suppose that $u$ and $w$ are
basic at every order. Then $v$ is basic at every order.
\end{subsub}

\begin{proof} Let $N$ be a positive integer such that \m{N\gg 0} and \ \m{{\bf
N}=(N,\ldots,N)}. Suppose that \Nligne
\m{w\equiv(Q_1(\pi_1),\ldots,Q_n(\pi_n))\quad
(\text{mod.} \quad (\pi^N))}, where \m{Q_1,\ldots,Q_n\in\C[X]}. Let
\m{{\bf m}=(m_1,\ldots,m_n)} \ be an $n$-tuple of positive integers, and \
\m{v=(v_i)_{1\leq i\leq n}}. Suppose that
\[u \ \equiv \ (P_1(\pi_1),\ldots,P_n(\pi_n)) \quad\quad (\text{mod.} \quad
\BS{\Pi}^{\bf N}) \]
Then we have
\[Q_i(\pi_i) \ \equiv \ P_i(\pi_i).v_i \quad\quad (\text{mod.} \quad
(\pi_i^N))\]
for \m{1\leq i\leq n}. We can write \m{P_i(X)} as \ \m{P_i(X)=X^{n_i}R_i(X)},
where \m{R_i(X)\in\C[X]} is such that \m{R_i(0)\not=0}. Then \m{Q_i(X)} is
also divisible by \m{X^{n_i}}: \m{Q_i(X)=X^{n_i}S_i(X)}, and we have in
\m{\ko_{\BS{\ka}x}} :
\[S_i(\pi_i) \ \equiv \ R_i(\pi_i).v_i  \quad\quad (\text{mod.} \quad
(\pi_i^{N'}))\]
for some integer \m{N'\gg 0}. We can write \ \m{R_i(X)=a_i.(1-X.T_i(X))},
with \ \m{a_i\in\C^*}, \m{T_i\in\C(X)}. We have then
\[v_i \ \equiv \
\frac{S_i(\pi_i)}{a_i}\sigg_{p=1}^{m_i-1}\big(\pi_iT_i(\pi_i)\big)^p
\quad\quad (\text{mod.} \quad \BS{\Pi}^{\bf m}) .\]
\end{proof}

\sepprop

For \m{1\leq i\leq n}, let \ \m{{\bf u}_{(i)}=((u_{(i)j})_{1\leq j\leq n}} \
be a generator of the ideal sheaf \m{\ki_{\kc_i}} of \m{\kc_i} in $\kc$, such
that for \m{1\leq j\leq n}, \m{u_{(i)j}\in\C[\pi_j]} (cf. corollary
\ref{ecl21}).

\sepprop

\begin{subsub}\label{ecl25}{\bf Proposition: } Let \m{v\in\ko_{\kc,x}}. then
$v$ is basic at every order if and only if for every $n$-tuple $\bf m$ of
positive integers, there exist an integer $q>0$ and
\m{P_1,\ldots,P_q\in\C[X]} such that
\[v \ \equiv \ \sigg_{1\leq j\leq q}P_j(\pi).{\bf u}^j_{(i)} \quad\quad
(\text{\em mod.} \quad \BS{\Pi}^{\bf m}) .\]
\end{subsub}

\begin{proof} We use the notations of the proof of lemma \ref{ecl23}.
Suppose that \m{v=(v_j)_{1\leq j\leq n}} is basic at every order. Let $N$ be a
positive integer and \ \m{{\bf N}=(N,\ldots,N)}. We will prove by induction on
$q\geq 0$ that we can write $v$ as
\begin{equation}\label{equ3}
v \ \equiv \ \sigg_{0\leq j\leq q}P_j(\pi).{\bf u}^j_{(i)} + \gamma_q{\bf
u}_{(i)}^{q+1} \quad\quad (\text{mod.} \quad \BS{\Pi}^{\bf N})
\end{equation}
with \m{P_0,\ldots,P_q\in\C[X]}, and \m{\gamma_q\in\ko_{\kc,x}}. This proves
proposition \ref{ecl25} if $q$ and $N$ are big enough.

For \m{q=0}, we have \ \m{v_i\equiv P(\pi_i) (\text{mod} \ \pi_i^N)},
for some \m{P\in\C[X]}, and we can take \m{P_0=P}. Suppose that the
result is true for $q$ and that we have \eqref{equ3}.
Since \ \m{v-\sigg_{1\leq j\leq q}P_j(\pi).{\bf u}^j_{(i)}} \ is basic at any
order, using the same method as in the proof of lemma \ref{ecl23}, we see that
\m{\gamma_q} is basic at order $\bf N'$, where \m{{\bf N'}=(N',\ldots,N')},
for some integer \m{N'\gg 0}. As in the case \m{q=0} we have
\[\gamma_q \ \equiv \ P_{q+1}(\pi)+{\bf u}_{(i)}.\gamma_{q+1} \quad\quad
(\text{mod} \quad \BS{\Pi}^{\bf N'}) ,\]
with \m{P_{q+1}\in\C[X]}. Hence
\[v \ \equiv \ \sigg_{0\leq j\leq q+1}P_j(\pi).{\bf u}^j_{(i)} +
\gamma_{q+1}{\bf u}_{(i)}^{q+2} \quad\quad (\text{mod.} \quad \BS{\Pi}^{\bf
N})\]
\end{proof}

\sepprop

\begin{subsub}\label{ecl24}{\bf Proposition: } Let \
\m{\alpha=(\alpha_1,\ldots,\alpha_n)\in\ko_{\kc,x}} \ be such that there
exists\Nligne \m{P_1,\ldots,P_{n-1}\in\C[X]} such that, for \m{1\leq
i\leq n-1}, we have \ \m{\alpha_i\equiv P_i(\pi_i)} (mod. \m{(\pi_i^{q_i})}).
Then there exists \m{P_n\in\C[X]} such that \ \m{\alpha_n\equiv
P_n(\pi_n)} (mod. \m{(\pi^{q_n})}), i.e. $\alpha$ is a basic element of
order $\bf q$.
\end{subsub}

\begin{proof} By induction on $n$. The case \m{n=2} is an easy consequence of
proposition \ref{ecl4}. Suppose that \m{n\geq 3} and that the result is true
for \m{n-1}.

By subtracting multiples of \m{(0,\ldots,0,\pi_i^{q_i},0,\ldots,0)} we may
assume that for \Nligne\m{1\leq i\leq n-1}, \m{\alpha_i\in\C[\pi_i]}. By
subtracting  a regular function on a neighborhood of $C$ in $\kc$, and a
multiple of \m{(\pi_1^{q_1},0,\ldots,0)} we may also assume that
\m{\alpha_1=0}. The ideal sheaf of
\m{\kc_1} is generated by \m{{\bf u}_{(1)}}. We can then write \
\m{\alpha=\beta{\bf u}_{(1)}}, with \ \m{\beta=(\beta_i)_{1\leq
i\leq n}\in\ko_{\kc,x}}. We have
\[(\alpha_2,\ldots,\alpha_{n-1})=(\beta_2,\ldots,\beta_{n-1}).
(u_{(1)2},\ldots,u_{(1)n-1}) \ , \]
hence by lemma \ref{ecl23}, \m{(\beta_2,\ldots,\beta_{n-1})} is a basic
element at any order. By the induction hypothesis, there exists
\m{Q\in\C[X]} such that \ \m{\beta_n\equiv Q(\pi_n)} (mod.
\m{(\pi_n^{q_n-p_{1n}})}). Since
\m{u_{(1)n}} is a multiple of \m{\pi_n^{p_{1n}}} (from the definition of
\m{p_{1n}}), it follows that \ \m{\alpha_n\equiv u_{(i)n}Q(\pi_n)} (mod.
\m{(\pi^{q_n})}).
\end{proof}

\end{sub}

\sepsub

\Ssect{Simple primitive curves and fragmented deformations}{ecl_sim}

Let \m{C_n} be a primitive multiple curve of multiplicity $n$ and associated
smooth curve $C$. Let \m{\ki_C} be the ideal sheaf of $C$ in \m{C_n}. It is
obvious from proposition \ref{ecl12}, 1-, that if there exists a fragmented
deformation of \m{C_n}, then we have \m{\ki_{C,C_n}\simeq\ko_{C_{n-1}}}, i.e.
\m{C_n} is {\em simple} (cf. \ref{cmprq}). Conversely we have

\sepprop

\begin{subsub}\label{theo_sim}{\bf Theorem:} Let \m{C_n} be a simple primitive
multiple curve of multiplicity $n$. Then there exists a fragmented deformation
of \m{C_n}.
\end{subsub}

\begin{proof} According to theorem \ref{theo_121}, there exists a flat family
of smooth projective curves \ \m{\tau:\kc\to \C} \ such that
\m{\tau^{-1}(0)\simeq C} and that \m{C_n} is isomorphic to the $n$-th
infinitesimal neighborhood of $C$ in $\kc$. Let \m{\rho_n:\C\to\C} be the map
defined by \m{\rho_n(z)=z^n}, and \ \m{\theta=\rho_n\circ\tau:\kc\to\C}. It is
a flat morphism, \m{\theta^{-1}(0)=C_n}, and for every \m{z\not=0} in the 
image of $\tau$, \m{\theta^{-1}(z)} is a disjoint union of $n$ smooth 
irreducible curves. We can then apply the process of proposition 
\ref{compconn3} to obtain the desired fragmented deformation: it is \ 
\m{\kc\times_\C\C}
\xmat{\kc\times_\C\C\ar[r]^-\pi\ar[d] & \C\ar[d]^{\rho_n} \\
\kc\ar[r]^-\theta & \C }
\end{proof}

\sepprop

\begin{subsub}\label{rem_sim}{\bf Remark: }\rm let \m{(p_{ij})} be the spectrum
of the fragmented deformation constructed in the proof of theorem
\ref{theo_sim}. Then it is easy to see that \m{p_{ij}=1} for \m{1\leq i,j\leq
n}, \m{i\not=j}. If \m{x\in C}, then \
\m{(\kc\times_\C\C)_x=\ko_{\kc,x}\ot_{\ko_{\C,x}}\ko_{\C,x}}, and if
\m{t=I_\C\in\ko_{\C,x}}, we have for \m{1\leq k\leq n}
\[(\pi_1,\ldots,\pi_{k-1},0,\pi_{k+1},\ldots,\pi_n) \ = \
\frac{1}{n-1}(1\ot t-e^{\frac{2ki\pi}{n}}(t\ot 1)) \ . \]
\end{subsub}

\end{sub}

\sepsec

\section{Stars of a curve}\label{stars}

\Ssect{Definitions}{s-Def}

Let $S$ be a smooth irreducible curve, and \m{P\in S} (we can also take for
\m{(S,P)} the germ of a smooth curve). Let $n$ be a positive integer.

\sepprop

\begin{subsub}\label{s_def}{\bf Definition:} An {\em $n$-star} (or more simply,
a {\em star}) of \m{(S,P)} is an algebraic variety \m{\BS{\ks}} such that
\begin{enumerate}
\item[(i)] $\BS{\ks}$ is the union of $n$ irreducible components
$S_1,\ldots,S_n$, with fixed isomorphisms $S_i\simeq S$, $1\leq i\leq n$.
\item[(ii)] For $1\leq i<j\leq n$, $S_i\cap S_j$ has only one closed point,
namely $P$.
\item[(iii)] There exists a morphism \ $\pi:\BS{\ks}\to S$, such that for 
$1\leq i\leq n$, the restriction \Nligne $\pi_{\mid S_i}:S_i\to S$ \ is the 
isomorphism $S_i\simeq S$ of (i).
\end{enumerate}
\end{subsub}
All the $n$-stars of \m{(S,P)} have the same underlying Zariski topological
space \m{S(n)} and set of closed points. The latter is \ \m{(\bigcup_{1\leq
i\leq n}\widehat{S}_i)/\sim}, where \m{\widehat{S}_i} is the set of closed
points of \m{S_i}, and the equivalence relation $\sim$ is defined by: for
\m{x\in\widehat{S_i}} and \m{y\in\widehat{S_j}}, \m{x\sim y} if and only if 
\m{i=j} and \m{x=y}, or \m{x=P\in\widehat{S_i}} and \m{y=P\in\widehat{S_j}}. 
An open subset of \m{\BS{\ks}} is defined by open subsets \m{U_1} of 
\m{S_1},$\ldots$, \m{U_n} of \m{S_n}, such that for
\m{1\leq i<j\leq n}, we have \m{P\in U_i} if and only if \m{P\in U_j}.

\sepprop

The {\em initial} star \m{\BS{\ks}_0} of \m{(S,P)} is defined as follows: for
every open subset $U$ of \m{S(n)}, \m{\ko_{\BS{\ks}_0}(U)} is the set of \
\m{(\alpha_1,\ldots,\alpha_n)\in\ko_{S_1}(U\cap
S_1)\times\cdots\ko_{S_n}(U\cap S_n)} \ such that if \m{P\in U} then
\m{\alpha_1(P)=\cdots=\alpha_n(P)} .

For every $n$-star \m{\BS{\ks}} of \m{(S,P)}, there is a unique dominant
morphism \ \m{\BS{\ks}_0\to \BS{\ks}} inducing the identity on each component.
So \m{\ko_{\BS{\ks},P}} is a subring of \m{\ko_{\BS{\ks}_0,P}}.

Note that (iii) is equivalent to
\begin{enumerate}
\item[(iii)'] For every $\alpha\in\ko_{S,P}$, we have \ $(\alpha,\ldots,\alpha)
\in\ko_{\BS{\ks},P}$.
\end{enumerate}

\sepprop

\begin{subsub}\label{s_def2}{\bf Definition:} An {\em oblate $n$-star} (or
more simply, an {\em oblate star}) of \m{(S,P)} is an $n$-star
\m{\BS{\ks}} such that some neighborhood of $P$ in \m{\BS{\ks}} can be
embedded in a smooth surface.
\end{subsub}

\sepprop

\begin{subsub}\label{s_prop1}{\bf Proposition:} An $n$-star \m{\BS{\ks}} is
oblate if and only if \ \m{\pi^{-1}(P)\simeq\spec(\C[X]/(X^n))}.
\end{subsub}

(cf. prop. \ref{ecl2b}).

\sepprop

Let \ \m{I\subset\{1,\ldots,n\}} \ be a nonempty subset. Let \
\m{\BS{\ks}^{(I)}=\bigcup_{i\in I}S_i\subset \BS{\ks}}. If \m{\BS{\ks}} is
oblate then
\m{\BS{\ks}^{(I)}} is oblate too.

\end{sub}

\sepsub

\Ssect{Properties of oblate stars}{P-star}

Let \m{\BS{\ks}} be an oblate $n$-star of $S$.
Recall that $t$ denotes a generator of the maximal ideal of $P$ in $S$. We
will denote this generator on \m{S_i\subset \BS{\ks}} by \m{t_i}. We will also
denote by $\pi$ the element \m{t\circ\pi} of the maximal ideal of $P$ in
\m{\BS{\ks}}. Let \m{\ki_P} be the ideal sheaf of $P$ in \m{\BS{\ks}}.

We begin with 2-stars:

\sepprop

\begin{subsub}\label{s_prop2} {\bf Proposition:} Suppose that \m{n=2}. Then

{\bf 1 -- } There exists a unique integer \m{p>0} such that \m{\ki_{P,P}/(\pi)}
is generated by the image of \m{(t_1^p,0)}.

{\bf 2 -- } The image of \m{(0,t_2^p)} is also a generator of
\m{\ki_{P,P}/(\pi)}.

{\bf 3 -- } \m{(0,t_2^p)} (resp. \m{(t_1^p,0)}) is a generator of the ideal
sheaf of \m{S_1} (resp. \m{S_2}) at $P$.

{\bf 4 -- } \m{\ko_{S^{(2)},P}} consists of pairs \
\m{(\alpha,\beta)\in\ko_{S,P}\times\ko_{S,P}} \ such that \
\m{\alpha-\beta\in(t^p)}.
\end{subsub}

(cf. prop. \ref{ecl3} and \ref{ecl4}).

\sepprop

Now suppose that \m{n\geq 2}. Let \ \m{I=\{i,j\}\subset\{1,\ldots,n\}}, with
\m{i\not=j}. Then \m{S_i\cup S_j\subset \BS{\ks}} is a $2$-star of $S$. Hence
by proposition \ref{s_prop2} there exists a unique integer \m{p_{ij}>0} such
that \m{\ki_{P,P}/(\pi)} (on \m{S_i\cup S_j}) is generated by the image of
\m{(t_i^{p_{ij}},0)} (and also by the image of \m{(0,t_j^{p_{ij}})}). Let 
\m{p_{ii}=0}. Then the symmetric matrix
\m{(p_{ij})_{1\leq i,j\leq n}} is called the {\em spectrum} of \m{\BS{\ks}}.

There exists an element \m{v_{ij}=(\nu_m)_{1\leq m\leq n}} such that
\m{\nu_i=0} and \m{\nu_j=t_j^{p_{ij}}}. For every integer $m$ such that \
\m{1\leq m\leq n}, \m{m\not=i,j}, there exists an invertible element
\m{\beta_{ij}^{(m)}\in\ko_{S,P}} such that \
\m{\nu_m=\beta_{ij}^{(m)}t_m^{p_{im}}}. Let \m{\beta_{ij}^{(i)}=0},
\m{\beta_{ij}^{(j)}=1}.

\sepprop

\begin{subsub}\label{s_prop3} {\bf Proposition:}
Let \ \m{{\bf b}^{(m)}_{ij}=\beta_{ij}^{(m)}(P)\in\C}. Then we have, for all
integers \m{i,j,k,m,q} such that \m{1\leq i,j,k,m,q\leq n}, \m{i\not=j},
\m{i\not=k}
\[{\bf b}_{ik}^{(m)}{\bf b}_{ij}^{(q)} \ = \ {\bf b}_{ik}^{(q)}{\bf
b}_{ij}^{(m)} .\]
In particular we have \ \m{{\bf b}_{ij}^{(m)}={\bf b}_{ik}^{(m)}{\bf
b}_{ij}^{(k)}} \ and \ \m{{\bf b}_{ij}^{(m)}{\bf b}_{im}^{(j)}=1}.

For all distinct integers
$i,j,k$ such that \m{1\leq i,j,k\leq n}, we have
\[{\bf b}_{ki}^{(j)} \ = \ - {\bf b}_{ik}^{(j)}{\bf b}_{ji}^{(k)} \ . \]
\end{subsub}

(cf. prop. 4.3.2 and 4.4.6).

\sepprop

Let $p$ be an integer such that \m{1\leq p<n}, and 
\m{(i_1,j_1),\ldots,(i_p,j_p)}
$p$ pairs of distinct integers of \m{\{1,\ldots,n\}}. Then the image of \
\m{\prod_{m=1}^p{\bf v}_{i_mj_m}} \ is a generator of \
\m{(\ki_{P,P}^p+(\pi))/(\ki_{P,P}^{p+1}+(\pi))}.

Let \ \m{I\subset\{1,\cdots,n\}} \ be a nonempty subset, distinct from
\m{\{1,\cdots,n\}}. Let \ \m{i\in\{1,\cdots,n\}\backslash I}. Let
\[{\bf v}_{I,i} \ = \ \prod_{j\in I} {\bf v}_{ji} \ .\]

\sepprop

\begin{subsub}\label{s_theo1}{\bf Proposition: } The ideal sheaf of
\m{\BS{\ks}^{(I)}} in \m{\BS{\ks}} is generated by \m{{\bf v}_{I,i}} at $P$.
\end{subsub}

(cf. prop. 4.3.3).

\sepprop

Note that if \ \m{I=\{1,\cdots,n\}\backslash\{i\}} \ then \ \m{{\bf
v}_{I,i\mid S_j}=0} \ if \m{j\not=i}, and \ \m{{\bf v}_{I,i\mid
S_i}=t_i^{q_i}}, with \ \m{q_i=\sigg_{1\leq j\leq n}p_{ij}}.

Let $i$ be an integer such that \m{1\leq i\leq n} and
\m{J_i=\{1,\ldots,n\}\backslash\{i\}}.
Let \m{\kk_i} be the image of \m{\ko_{\BS{\ks}}} in \m{\prod_{1\leq j\leq
n,j\not=i}\ko_{S_j}/(t_j^{q_j})}. We can view \m{\kk_i} as a $\C$-algebra.
For every
\m{\alpha=(\alpha_m)\in\ko_{\BS{\ks},P}}, let \m{k_i(\alpha)} be the image of
$\alpha$ in \m{\kk_i}.

\sepprop

\begin{subsub}\label{s_prop4}{\bf Proposition: } There exists a morphism of
$\C$-algebras
\[\Psi_i : \kk_i\lra\ko_{S_{i,P}}/(t_i^{q_i})\]
such that for every \m{(\alpha_m)_{1\leq m\leq
n,m\not=i}\in\ko_{\BS{\ks}^{(J_i)},P}}, \m{\alpha_i\in\ko_{S_{j,P}}}, we
have \m{\alpha=(\alpha_m)_{1\leq m\leq n}\in\ko_{\BS{\ks},P}} if and only if \
\m{\Psi_i(k_i(\alpha))=[\alpha_i]_{q_i}}.
\end{subsub}

(cf. prop. 4.4.1).

\sepprop

The morphism \m{\Psi_i} has the following properties:
\begin{enumerate}
\item[(i)] For every \m{(\alpha_m)_{1\leq m\leq
n,m\not=i}\in\ko_{\BS{\ks}^{(J_i)},P}}, we have 
$\Psi_i(\alpha)(P)=\alpha_m(P)$ for $1\leq m\leq n, m\not=i$.
\item[(ii)] We have \ $\Psi_i((t_m)_{1\leq m\leq n,m\not=i})=t_i$.
\item[(iii)] Let \ $j,k\in\{1,\cdots,n\}$ \ be such that $i,j,k$ are distinct.
Let $\bf w$ be the image of ${\bf v}_{jk}$ in $\kb_i$.
Then there exists $\lambda\in\ko_{S_{i,P}}^*$ such that \
$\Psi_i({\bf w})=\lambda t_i^{p_{ij}}$.
\item[(iv)] Let $j$ be an integer such that $1\leq j\leq n$ and $j\not=i$. Let
$\bf w$ be the image of ${\bf v}_{ij}$ in $\kk_ix$. Then we have \
$\ker(\Psi_i)=({\bf w})$.
\end{enumerate}

\sepprop

\begin{subsub}\label{s_conv} Converse -- \rm Let $\BS{\ks}^{[n-1]}$ be a
\m{(n-1)}-star of $S$, with components  \m{S_1,\ldots,S_{n-1}}, of spectrum
\m{(p_{jk})_{1\leq j,k\leq n-1}}. Let \m{p_{nj}=p_{jn}}, \m{1\leq j<n} be 
positive integers, and \m{p_{nn}=0}. For \m{1\leq j\leq n}, let \ 
\m{q_j=\sigg_{1\leq k\leq n}p_{kj}}.

Let \m{S_n} be another copy of $S$. Let \m{\kk_n} be the image of
\m{\ko_{\BS{\ks}^{[n-1]}}} in \ \m{\prod_{1\leq j\leq
n-1}\ko_{S_j}/(t_j^{q_j})} \
and
\[\Psi_n : \kk_n\lra\ko_{S_n}/(t_n^{q_n})\]
a morphism of $\C$-algebras satisfying properties
(i), (ii), (iii) above. Let $\kk$ be the subsheaf of algebras of
\m{\ko_{\BS{\ks}_0}} defined by: \m{\kk=\ko_{\BS{\ks}_0}} on
\m{\BS{\ks}_0\backslash\{P\}}, and for every \ \m{\alpha=(\alpha_m)_{1\leq
m\leq
n}\in\ko_{\BS{\ks}_0,P}}, \m{\alpha\in\kk_P} if and only if \
\m{\Psi_n(\alpha')=\lbrack\alpha_n\rbrack_{q_n}} (where $\alpha'$ is the image
of \m{(\alpha_m)_{1\leq m\leq n-1}} in \m{\kk_n}).

It is easy to see that $\kk$ is the structural sheaf of an oblate $n$--star of
$S$.
\end{subsub}

\sepprop

Let \ \m{\kh=\prod_{1\leq j\leq n}(t_j^{q_j-1})/(t_j^{q_j})\simeq\C^n} \
and \m{\kk} be the image of \m{\ko_{\BS{\ks}}} in \m{\prod_{1\leq j\leq
n}\ko_{S_j}/(t_j^{q_j})}. We can view \m{\kk} as a $\C$-algebra. Let \
\m{\kj=\kh\cap\kk}.

\sepprop

\begin{subsub}\label{s_prop5}{\bf Proposition: } There exists a unique \
\m{\lambda(\BS{\ks})=(\lambda_1,\ldots,\lambda_n)\in\P_n(\C)} such that for
every \m{{\bf u}=(u_j)_{1\leq j\leq n}\in\kh}, we have \m{{\bf
u}\in\kj} if and only if \ \m{\lambda_1u_1+\cdots\lambda_nu_n=0}. The
\m{\lambda_i} are all non zero.
\end{subsub}

(cf. prop. 4.4.5).

\sepprop

For all distinct integers $i$, $j$ such that \m{1\leq i,j\leq n}, we have
\[
\frac{\lambda_i}{\lambda_j} \ = \ - \prod_{1\leq m\leq
n,m\not=i,j}{\bf b}^{(j)}_{mi} .
\]
\end{sub}

\sepsub

\Ssect{Construction of oblate stars of a curve}{const-st}

Consider an oblate \m{(n-1)}-star of $S$, \m{\BS{\ks}^{[n-1]}}, with \m{n-1}
irreducible components \m{S_1,\ldots,S_{n-1}}, copies of $S$. Let
\m{(p_{ij}^{[n-1]})_{1\leq i,j<n}} be its spectrum. For \m{1\leq i<n}, let \
\m{q_i^{[n-1]}=\sigg_{1\leq j<n}p_{ij}^{[n-1]}}. We denote by
\m{\ki_P^{[n-1]}} the ideal of $P$ in \m{\ko_{\BS{\ks}^{[n-1]},P}}. Let \
\m{\lambda(\BS{\ks}^{[n-1]})=(\lambda_1,\ldots,\lambda_{n-1})}.

Let \m{p_{1n},\ldots,p_{n-1,n}} be positive integers,
\m{q_i=q_i^{[n-1]}+p_{in}} \ for \m{1\leq i<n}, and \Nligne
\m{q_n=p_{i1}+\cdots+p_{n-1,n}}. Let \m{{\bf u}\in\ki_{P,P}^{[n-1]}} whose
image generates \m{\ki_P^{[n-1]}/((\ki_P^{[n-1]})^2+(\pi))}, of the form
\[{\bf u} \ = \ (\beta_1t_1^{p_{1n}},\ldots,\beta_{n-1}
t_{n-1}^{p_{n-1,n}}) ,\]
with \m{\beta_i\in\ko_{S_{i,P}}} invertible for \m{1\leq i<n}.

Let \m{\kk^{[n-1]}} be the image of \m{\ko_{\BS{\ks}^{[n-1]}}} in \
\m{\ko_{S_1}/
(t_1^{q_1})\times\cdots\times\ko_{S_{n-1}}/ (t_{n-1}^{q_{n-1}})} . We will
also denote by $\bf u$ the image of $\bf u$ in \m{\kk^{[n-1]}}. Let \
\m{\kq=\kk^{[n-1]}/({\bf u})} , \m{\rho:\kk^{[n-1]}\to\kq} \ the projection and
\ \m{t_n=\rho(\pi)}.

\sepprop

\begin{subsub}{\bf Proposition: }\label{s_prop6}{\bf 1 --} We have \
\m{t_n^{q_n}=0} .

{\bf 2 --} We have \
\m{t_n^{q_n-1}=0} \ if and only if
\[\frac{\lambda_1}{\beta_1(P)}+\cdots+
\frac{\lambda_{n-1}}{\beta_{n-1}(P)} \ = \ 0 \ .\]

We suppose now that \ \m{\frac{\lambda_1}{\beta_1(P)}+\cdots+
\frac{\lambda_{n-1}}{\beta_{n-1}(P)}\not=0}. Then

{\bf 3 --} For every \ \m{\epsilon\in\kk^{[n-1]}} \ such that \
\m{\epsilon(P)\not=0}, we have \ \m{t^{q_n-1}\epsilon\not\in({\bf u})}.

{\bf 4 --} For every \ \m{\eta\in\kk^{[n-1]}/({\bf u)}}, and every integer
$k$ such that \m{1\leq k<q_n}, we have \ \m{t_n^k\eta=0} \ if and only if
$\eta$ is a multiple of \m{t_n^{q_n-k}}.

{\bf 5 --} \m{\kk^{[n-1]}/({\bf u)}} is a flat
\m{\C[t_n]/(t_n^{q_n})}-module.
\end{subsub}

(cf. prop. \ref{ecl17}).

\sepprop

\begin{subsub}\label{cst_st} Construction of stars of a curve -- \rm
Suppose that \ \m{\frac{\lambda_1}{\beta_1(P)}+\cdots+
\frac{\lambda_{n-1}}{\beta_{n-1}(P)}\not=0}. From proposition \ref{s_prop6},
5-, it is easy to prove, using \ref{s_conv}, that there is a unique oblate
$n$-star \m{\BS{\ks}} such that \m{\BS{\ks}^{[n-1]}} is the union \
\m{\bigcup_{1\leq
i\leq n-1}S_i} \ in \m{\BS{\ks}} and \m{\Psi_n} is the quotient map \
\m{\kk_n=\kk^{[n-1]}\to\kq}.
\end{subsub}
\end{sub}

\sepsub

\Ssect{Morphisms of stars}{s-morph}

Recall that if $\BS{\ks}$ is an oblate $n$-star of $S$, then we have a
canonical inclusion of sheaves of algebras (on the underlying topological
space \m{S(n)} of \m{\BS{\ks}}) \ \m{\ko_{\BS{\ks}}\subset\ko_{\BS{\ks}_0}}.

Let $\BS{\ks}$, \m{\BS{\ks}'} be oblate $n$-stars of $S$, with irreducible
components \m{S_1,\ldots,S_n}, and \m{f:\BS{\ks}\to\BS{\ks}'} a morphism
inducing the identity on all the components. Such a morphism exists if and
only if \ \m{\ko_\BS{\ks'}\subset\ko_\BS{\ks}}, and in this case $f$ is unique 
and is induced by the previous inclusion. Let \m{(p_{ij})} (resp. 
\m{(p_{ij}'))} be the spectrum of $\BS{\ks}$ (resp. \m{\BS{\ks}'}).

\sepprop

\begin{subsub}\label{s_prop7}{\bf Proposition: } We have \ \m{p_{ij}\leq
p'_{ij}} \ for \m{1\leq i,j\leq n}. If $f$ is not the identity morphism then
there exist $i$, $j$ such that \ \m{p_{ij}<p'_{ij}}.
\end{subsub}

\begin{proof}
Let \m{I=\{i,j\}}. Then $f$ induces a morphism
\m{\BS{\ks}^{(I)}\to{\BS{\ks}'}^{(I)}}. So we have \
\m{\ko_{{\BS{\ks}'}^{(I)},P}\subset\ko_{\BS{\ks}^{(I)},P}}. From proposition
\ref{s_prop2}, 4-, it follows that \m{p_{ij}\leq p'_{ij}}.

Suppose now that \m{p'_{ij}=p_{ij}} for \m{1\leq i,j\leq n}. We must prove
that \m{\BS{\ks}=\BS{\ks}'}, i.e. that \
\m{\ko_{{\BS{\ks}'},P}=\ko_{\BS{\ks},P}}. This is done by
induction on $n$. For \m{n=2} it is obvious. Suppose that it is true for
\m{n-1}. Let \m{I=\{1,\ldots,n-1\}}. Then $f$ induces a morphism \
\m{f_{n-1}:\BS{\ks}^{(I)}\to{\BS{\ks}'}^{(I)}}. It follows from the induction
hypothesis that \ \m{\BS{\ks}^{(I)}={\BS{\ks}'}^{(I)}}. Since the integers
\m{q_i} are the same for $\BS{\ks}$ and \m{\BS{\ks}'}, the algebras \m{\kk_n}
for $\BS{\ks}$ and \m{\BS{\ks}'} (cf. proposition \ref{s_prop4}) are also the
same. Now let \ \m{\alpha\in\ko_{\BS{\ks},P}}, and let \m{\beta\in\kk_n}
be the image of $\alpha$. Let \m{\alpha'\in\ko_{\BS{\ks}',P}} be such
that its image in \m{\kk_n} is also $\beta$. Then \m{\alpha-\alpha'} belongs
to the ideal generated by the \m{(0,\ldots,0,t_i^{q_i},0\ldots,0)}, \m{1\leq
i\leq n}, which is included in \m{\ko_{\BS{\ks}',P}}. Hence \
\m{\alpha\in\ko_{\BS{\ks}',P}}.
\end{proof}

\sepprop

\begin{subsub}\label{s-lem1}{\bf Lemma: } Suppose that $f$ is not the identity
morphism. Then there exist an ideal \ \m{\ki\subset\ko_{\BS{\ks}',P}}
\ and \m{u\in\ki}, \m{v\in\ko_{\BS{\ks},P}} \ such that
\[u\ot v\not=0 \quad\quad \text{in} \quad\quad
\ki\ot_{\ko_{\BS{\ks}',P}}\ko_{\BS{\ks},P}\]
and \m{uv=0}.
\end{subsub}

\begin{proof}
Let \ \m{q_1=\sigg_{i=1}^np_{1i}}, \m{q'_1=\sigg_{i=1}^np'_{1i}}. According to
proposition \ref{s_prop7} we can assume that \ \m{q_1<q'_1}. Let \m{u}
be a generator of the ideal of \m{S_1} in \m{\ko_{\BS{\ks}',P}} and
\m{\ki=(u)}. Let \ \m{v=(t_1^{q_1},0,\ldots,0)}. We have \m{uv=0}. We have to
prove that \ \m{u\ot v\not=0}. We need only to find an
\m{\ko_{\BS{\ks}',P}}-module $M$ and a \m{\ko_{\BS{\ks}',P}}-bilinear map
\[\phi:\ki\ot_{\ko_{\BS{\ks}',P}}\ko_{\BS{\ks},P}\lra M\]
such that \ \m{\phi(u\ot v)\not=0}. We take \ \m{M=\ko_{S_1,P}/(t_1^{q'_1})},
which is a quotient of \m{\ko_{{\BS{\ks}'}}}. It is easy to verify that
\xmat{\phi:((\lambda_i)_{1\leq i\leq n}u,(w_i)_{1\leq i\leq n})\fmaps[r] &
\lambda_1w_1 \ (\text{mod} \ t_1^{q'_1})}
is well defined, bilinear, and that \ \m{\phi(u\ot v)\not=0}.
\end{proof}

\sepprop

\begin{subsub}\label{s-cor1}{\bf Corollary: } Suppose that $f$ is not the
identity morphism. Let $Y$ be an algebraic variety and \ \m{g:Y\to S} \ a
morphism such that \ \m{g^*:\ko_{\BS{\ks},P}\to\ko_{Y,P}} \ is injective. Then
\ \m{f\circ g:Y\to S'} \ is not flat.
\end{subsub}

\begin{proof}
We use the notations of the proof of lemma \ref{s-lem1}. We have a commutative
diagram
\xmat{\ko_{\BS{\ks},P}\ar[rr]^-{g^*}\ar[d]^{\lambda_S} & &
\ko_{Y,P}\ar[d]^{\lambda_Y} \\
\ki\ot_{\ko_{\BS{\ks}',P}}\ko_{\BS{\ks},P}\ar[rr]^-{I_\ki\ot g^*}
\ar[d]^{\mu_S} & & \ki\ot_{\ko_{\BS{\ks}',P}}\ko_{Y,P}\ar[d]^{\mu_Y} \\
\ko_{\BS{\ks},P}\ar[rr]^-{g^*} & & \ko_{Y,P}}
where \m{\lambda_S(\alpha)=u\ot\alpha}, \m{\mu_S(u\ot\alpha)=u\alpha}, and
\m{\lambda_Y,\mu_Y} are defined similarly. It follows that \ \m{\mu_Y(u\ot
g^*v)=0}. We will show that \ \m{u\ot g^*v\not=0}, and this will imply that
\m{f\circ g} is not flat. Let \ \m{w=(t_1^{q'_1},0,\ldots,0)}. Then we have \
\m{\ki\simeq\ko_{\BS{\ks}',P}/(w)}, and from the exact sequence of
\m{\ko_{\BS{\ks}',P}}-modules \ \m{0\to(w)\to\ko_{\BS{\ks}',P}\to\ki\to 0} we
deduce that \ \m{\ker(\lambda_Y)=(w).\ko_{Y,P}}. Suppose that \m{u\ot g^*v=0}.
Then \m{g^*v} is a multiple of $w$: \m{g^*v=w.a}, for some \m{a\in\ko_{Y,P}}.
But we have \m{w=g^*\pi^{q'_1-q_1}v}. Hence \
\m{g^*v.(1-g^*\pi^{q'_1-q_1})=0}.
Since \m{1-g^*\pi^{q'_1-q_1}} is invertible, we have \m{g^*v=0}, which is
false since \m{g^*} is injective. Hence \ \m{u\ot g^*v\not=0}.
\end{proof}

\end{sub}

\sepsub

\Ssect{Structure of ideals}{str_ide}

Let $\BS{\ks}$ be an oblate $n$-star of $S$.

\sepprop

\begin{subsub}{\bf Proposition: }\label{s-prop8} Let
\m{\ki\subset\ko_{\BS{\ks},P}} be a proper ideal. Then

{\bf 1 - } There exists a positive integer $k$ such that \m{k\leq n} and a
filtration by ideals
\[\nsp=\ki_{k+1}\subset\ki_k\subset\cdots\subset\ki_1=\ki\]
such that, for \m{1\leq i\leq k} there exists a positive integer $j$ such that
\m{j\leq n} and an isomorphism \ \m{\ki_i/\ki_{i+1}\simeq\ko_{S_j,P}} of
\m{\ko_{\BS{\ks},P}}-modules.

{\bf 2 - } If \ \m{\ki_i/\ki_{i+1}\simeq\ko_{S_j,P}}, then \ \m{\ki_{i+1}
\subset\ki_{S_j}} \ and \ \m{\ki_i\not\subset\ki_{S_j}}.
\end{subsub}

\begin{proof} We prove {\bf 1-} by induction on $n$. The case \m{n=1} is
trivial. Suppose that \m{n>1} and that the result is true for \m{n-1}. Let
\m{\kj_1} be the ideal sheaf of \m{S_1\subset\BS{\ks}}, and \
\m{\BS{\ks}'=S_2\cup\cdots\cup S_{n-1}\subset\BS{\ks}}. We can view \m{\kj_1}
as an ideal of \m{\ko_{\BS{\ks}',P}}. We can suppose that
\m{\ki\not\subset\ko_{\BS{\ks}',P}}, i.e that some element of $\ki$ has a
nonzero first coordinate. Let $m$ be the smallest positive integer such that
$\ki$ contains an element $u$ of the form
\[u \ = \ (t^m,\alpha_2,\ldots,\alpha_n) \ .\]
Then every element $v$ of $\ki$ can be written as
\[v \ = \ \lambda u+v' \ , \]
with \ \m{\lambda\in\ko_{\BS{\ks},P}} \ and \ \m{v'\in\kj_1\cap\ki}, and the
first coordinate of $\lambda$ is uniquely determined. It follows that \
\m{\ki/(\kj_1\cap\ki)\simeq\ko_{S_1,P}}. We can apply the recurrence
hypothesis to the ideal \m{\kj_1\cap\ki} of \m{\ko_{\BS{\ks}',P}} and get a
filtration of it, from which we deduce the filtration of $\ki$. This proves
{\bf 1-} for $n$.

Now we prove {\bf 2-}. Let \ \m{\alpha\in\ko_{\BS{\ks},P}\backslash\ki_{S_j}}.
Let \m{u\in\ki_i} be over a generator of \m{\ki_i/\ki_{i+1}}. Then the image
of \m{\alpha u} in \m{\ki_i/\ki_{i+1}} is not zero, i.e. \ \m{\alpha
u\not\in\ki_{i+1}}. Hence \m{\alpha\not\in\ki_{i+1}}, and \ \m{\ki_{i+1}
\subset\ki_{S_j}}. Let
\m{v_i=(0,\ldots,0,t_i^{q_i},0,\ldots,0)\in\ko_{\BS{\ks},P}}. Then the image
of \m{v_i u} in \m{\ki_i/\ki_{i+1}} is not zero, hence \m{u\not\in\ki_{S_j}}
and \ \m{\ki_i\not\subset\ki_{S_j}}.
\end{proof}

\end{sub}

\sepsub

\Ssect{Star associated to a fragmented deformation}{st_fr_de}

We keep the notations of chapter \ref{Defecl}. Let \m{n\geq 2} be an integer,
\m{\pi:\kc\to S} a fragmented deformation of \m{C_n}, and
\m{\kc_1,\ldots,\kc_n} the irreducible components of $\kc$.

Recall that \m{S(n)} is the underlying (Zariski) topological space of any
$n$-star of $S$. Let \m{\kc^{top}} be the underlying topological space of
$\kc$. We have an obvious continuous map \ \m{\BS{\pi}:\kc^{top}\to S(n)}. Let
\m{\ka_n} be the sheaf of algebras on \m{S(n)} defined by: for every open
subset $U$ of \m{S(n)}, \m{\ka_n(U)} is the algebra of \
\m{(\alpha_1,\ldots,\alpha_n)\in\ko_\kc(\BS{\pi}^{-1}(U))} \ such that \
\m{\alpha_i\in\ko_{S_i}(U\cap S_i)} \ for \m{1\leq i\leq n}.

According to corollary \ref{ecl22}, for every \m{x\in C}, \m{\ka_{n,P}} is the
algebra of \m{(\alpha_1,\ldots,\alpha_n)\in\ko_{\kc,x}} such that
\m{\alpha_i\in\ko_{S,P}} for \m{1\leq i\leq n}.
\end{sub}

\sepprop

\begin{subsub}\label{st_fr_1}{\bf Proposition: } The sheaf \m{\ka_n} is the
structural sheaf of an oblate $n$-star of $S$.
\end{subsub}

\begin{proof} By induction on $n$. The case \m{n=1} is obvious. Suppose that
\m{n>1} and that the result is true for \m{n-1}. Let \
\m{\kc'=\kc_1\cup\cdots\kc_{n-1}\subset\kc}, and \m{\ka_{n-1}} the
corresponding oblate \m{(n-1)}-star of $S$. Let
\[\Phi_n:\kb_n\lra\ko_{\kc_n}/(\pi_n^{q_n})\]
be the morphism of proposition \ref{ecl9}. According to proposition
\ref{ecl24}, \m{\Phi_n} induces a morphism
\[\Psi_n:\kk_n\lra\ko_{S_n,P}/(t_n^{q_n}) \ .\]
By the definitions of \m{\ka_n} and \m{\Phi_n}, if \
\m{u=(\alpha_1,\ldots,\alpha_n)\in\ko_{S_1,P}\times\cdots\times\ko_{S_1,P}},
then \m{u\in\ka_{n,P}} if and only if \m{\Psi_n(u')=v}, where \m{u'} (resp.
$v$) is the image of $u$ in \m{\kk_n} (resp. \m{\ko_{S_n,P}/(t_n^{q_n})}). The
result follows then from \ref{s_conv}.
\end{proof}

\sepprop

We denote by \m{\ks(\kc)} (or more simply $\ks$) the oblate $n$-star
corresponding to $\ka_n$, so \m{\ko_{\ks(\kc)}=\ka_n}. From the definition of
\m{\ka_n} we get a canonical morphism
\[\BS{\Pi}:\kc\lra\ks\]
such that \ \m{\BS{\Pi}_{\mid\kc_i}=\pi_i:\kc_i\to S_i} \ for \m{1\leq i\leq
n}.

\sepprop

\begin{subsub}\label{st_fr_2}{\bf Theorem:} The morphism $\BS{\Pi}$ is flat.
\end{subsub}

\begin{proof} We need only to prove that $\BS{\Pi}$ is flat at any point $x$
of $C$. Let \m{\ki\subset\ko_{\ks,P}} be a proper ideal. We have to show that
the canonical morphism of \m{\ko_{\ks,P}}-modules
\[\tau=\tau_\ki:\ko_{\kc,x}\ot_{\ko_{\ks,P}}\ki\lra\ko_{\kc,x}\]
is injective. According to proposition \ref{s-prop8} there is a filtration by
ideals
\[\nsp=\ki_{k+1}\subset\ki_k\subset\cdots\subset\ki_1=\ki\]
such that, for \m{1\leq i\leq k} there exists a positive integer $j$ such that
\m{j\leq n} and an isomorphism \ \m{\ki_i/\ki_{i+1}\simeq\ko_{S_j,P}} of
\m{\ko_{\BS{\ks},P}}-modules. We will prove the injectivity of $\tau$ by
induction on $k$.

Recall that for \m{1\leq j\leq n}, \m{\ki_{S_j,P}=\ki_{S_j,\ks,P}} is a 
principal ideal, generated by an element \m{u_j} which is also a generator of \
\m{\ki_{\kc_j,x}=\ki_{\kc_j,\kc,x}} (cf. corollary \ref{ecl21} and proposition 
\ref{s_theo1}), and that the only zero coordinate of \m{u_j} is the $j$-th.

Suppose that \m{k=1}, so \m{\ki} is isomorphic to \m{\ko_{S_j,P}} for some
$j$. Let $u$ be a generator of $\ki$ and
\m{w\in\ko_{\kc,x}\ot_{\ko_{\ks,P}}\ki}, that can be written as \ \m{w=v\ot
u}, \m{v\in\ko_{\kc,x}}. Suppose that \ \m{\tau(v\ot u)=vu=0}. Since $\ki$ is
annihilated by \m{\ki_{S_j,P}}, we have \
\m{\ki\subset\big((0,\ldots,0,t_j^{q_j},0,\ldots,0)\big)}. Since \m{vu=0}, the
$j$-th component of $v$ is zero, i.e. \m{v\in\ki_{\kc_j,x}}. Hence $v$ is a
multiple of \m{u_j}: \m{v=\alpha u_j}. We have then
\begin{eqnarray*}
w & = & \alpha u_j\ot u\\
       & = & \alpha\ot u_ju \ \ \text{(because \ }u_j\in\ko_{\ks,P} )\\
       & = & 0 \ \ \ \ \text{(because \ }u_ju=0) \ .
\end{eqnarray*}
Hence $\tau$ is injective.

Suppose that the result is true for \m{k-1\geq 1} and that the filtration of
$\ki$ is of length $k$. According to proposition \ref{s-prop8}, {\bf 1-}, we
have \ \m{\ki/\ki_2\simeq\ko_{S_j,P}} for some $j$. Let \m{u\in\ki} be such
that its image in \m{\ki/\ki_2} is a generator, and \
\m{w\in\ko_{\kc,x}\ot_{\ko_{\ks,P}}\ki} \ such that \m{\tau(w)=0}. We can
write $w$ as \ \m{w=\alpha\ot v+\beta\ot u}, with
\m{\alpha,\beta\in\ko_{\kc,x}} and \m{v\in\ki_2}. Since \ \m{\alpha v+\beta
u=0}, we have \ \m{\beta u\in\ko_{\kc,x}\ki_2}, and \ \m{\ko_{\kc,x}\ki_2
\subset\ki_{\kc_j}} \ by proposition \ref{s-prop8}, {\bf 2-}, i.e. the $j$-th
coordinate of \m{\beta u} is zero. By proposition \ref{s-prop8}, {\bf 2-}, the
$j$-th coordinate of $u$ does not vanish, hence the $j$-th coordinate of
$\beta$ is zero, i.e. \m{\beta\in\ki_{\kc_j}}. Hence $\beta$ is a multiple of
\m{u_j}: \m{\beta=\gamma u_j}. We have then
\[\beta\ot u = \gamma u_j\ot u = \gamma\ot u_ju , \]
and \m{u_ju\in\ki_2} (because its image in \m{\ki/\ki_2} vanishes). It follows
that $w$ is the image of an element $w'$ of
\m{\ko_{\kc,x}\ot_{\ko_{\ks,P}}\ki_2}. We have \ \m{\tau_{\ki_2}(w')=0}, hence
by the induction hypothesis \m{w'=0}. It follows that we have also \m{w=0}.
\end{proof}

\sepprop

\begin{subsub}\label{s-rem}{\bf Remark: }\rm If \m{\ks'} is an oblate $n$-star 
of $S$, and if \ \m{\Pi':\kc\to\ks'} \ is a flat morphism compatible with the 
projections to $S$, then we have \m{\ks'=\ks(\kc)} and \m{\Pi'=\Pi}. This is 
an easy consequence of corollary \ref{s-cor1}.
\end{subsub}

\sepprop

\begin{subsub}\label{s-conv}Converse - \rm Let \m{\pi:\ks\to S} be an oblate
$n$-star of $S$. Let \m{\BS{\Pi}:\kc\to\ks} be a flat morphism such that for
every closed point \m{s\in\ks}, \m{\BS{\Pi}^{-1}(s)} is a smooth irreducible
projective curve. Let \m{C=\BS{\Pi}^{-1}(P)} and \
\m{\tau=\pi\circ\BS{\Pi}:\kc\to S}. Then \m{C_n=\tau^{-1}(P)} is a primitive
multiple curve of multiplicity $n$ and associated smooth curve $C$, and $\kc$
is a fragmented deformation of \m{C_n}. This is an easy consequence of
proposition \ref{ecl2b}.
\end{subsub}

\newpage

\section{Classification of fragmented deformations of length 2}\label{cl_fr_2}

Let \m{\pi:\kc\to\C} be a fragmented deformation of length 2. The
corresponding double curve \m{C_2} is \m{\pi^{-1}(0)}. Suppose that the
spectrum of $\kc$ is \m{\begin{pmatrix}0 & p \\ p & 0\end{pmatrix}}. This 
means that the infinitesimal neighborhoods of order $p$ of $C$ in \m{\kc_1} 
and \m{\kc_2} are isomorphic, i.e. we have an isomorphism of sheaves of 
algebras on $C$
\[\Phi : \ko_{\kc_1}/(\pi_1^p)\lra\ko_{\kc_2}/(\pi_2^p) \ , \]
and for every point $x$ of $C$, we have
\[\ko_{\kc,x} \ = \ \{(\alpha_1,\alpha_2)\in\ko_{\kc_1,x}\times\ko_{\kc_2,x} \
; \ \alpha_2 \ (\text{mod} \ \pi_2^p)=\Phi(\alpha_1 \ (\text{mod} \
\pi_1^p))\} \ .\]
Let \m{C_i^k} denote the infinitesimal neighborhood of order $k$ of $C$ in
\m{\kc_i}, \m{i=1,2}, \m{k>0}. It is a primitive multiple curve of
multiplicity $k$ and associated smooth curve $C$, and we have \m{C^p_1=C^p_2}.
Hence \m{C^{p+1}_1} and \m{C^{p+1}_2} appear as extensions of \m{C^p_1} in
primitive multiple curves of multiplicity \m{p+1}. According to \cite{dr1} and 
\cite{dr6} these extensions are classified by \m{H^1(C,T_C)} (\m{T_C} being 
the tangent sheaf on $C$). More precisely, we say that two such extensions 
$D$, \m{D'} are {\em isomorphic} if there exists an isomorphism \m{D\simeq D'} 
leaving \m{C^p_1} invariant. Then if $\kh$ is the set of isomorphism classes 
of such extensions, a bijection \ \m{\lambda:H^1(C,T_C)\to\kh} \ is defined in
\cite{dr1}, such that \m{\lambda(0)=C^{p+1}_1}.

On the other hand, it follows from \cite{ba_ei}, \cite{dr1} that the primitive
double curves with associated smooth curve $C$ and associated line bundle
\m{\ko_C} are classified by \ \m{\P(H^1(C,T_C))\cup\nsp} .

\sepprop

\begin{subsub}\label{theo_cf}{\bf Theorem: } The point of \
\m{\P(H^1(C,T_C))\cup\nsp} \ corresponding to \m{C_2} is
\m{\C.\lambda^{-1}(C_2^{p+1})}.
\end{subsub}

\begin{proof} According to \cite{dr1}, there exists an open covering
\m{(U_i)_{i\in I}} of $C$ such that for \m{k=1,2}, the open subset of
\m{C_k^{p+1}} corresponding to \m{U_i} is isomorphic to \
\m{U_i\times\spec(C[t]/(t^{p+1}))}. Here $t$ is \m{\pi_1} on $\kc_1$ and
\m{\pi_2} on $\kc_2$. We obtain then cocycles \m{(\theta_{ij}^{(k)})_{i,j\in
I}}, where \m{\theta_{ij}^{(k)}} is an automorphism of \
\m{U_{ij}\times\spec(C[t]/(t^{p+1}))}. We can also suppose that
\m{\omega_{C\mid U_i}} is trivial, for every \m{i\in I}. Let \m{dx_{ij}=dx} be
a generator of \m{\omega_C(U_{ij})}. Since the ideal sheaf of $C$ in
\m{C_k^{p+1}} is the trivial sheaf on \m{C_k^p}, we can write, using
the notations of \cite{dr1}, \m{\theta_{ij}^{(k)}=\phi_{\mu^{(k)}_{ij},1}} ,
with \ \m{\mu^{(k)}_{ij}\in\ko_C(U_{ij})[t]/(t^p)} , i.e. for every
\m{\alpha\in\ko_C(U_i)}, we have, at the level of regular functions
\[\theta_{ij}^{(k)}(\alpha) \ = \
\sigg_{m=0}^p\frac{1}{m!}(\mu^{(k)}_{ij}t)^m\frac{d^m\alpha}{dx^m} \ , \]
and \ \m{\theta_{ij}^{(k)}(t)=t} . Since \m{C_1^p=C_2^p} we can suppose that
\ \m{\mu_{ij}^{(1)}\equiv\mu_{ij}^{(2)} \ (\text{mod} \ t^{p-1})} . Hence \
\m{\tau_{ij}=\mu_{ij}^{(2)}-\mu_{ij}^{(1)}\in(t^{p-1})/(t^p)\simeq\ko_C(U_i)}.
The family \m{(\tau_{ij})} is (in some sense) a cocycle representing
\m{\lambda^{-1}(C_2^{p+1})} (cf. \cite{dr1}, \cite{dr6}).

We have \m{(\pi_1^{p+1})+(\pi_2^{p+1})\subset(\pi)} in \m{\ko_\kc}. Hence
\m{C_2=\pi^{-1}(0)} is contained in the subscheme $Z$ of $\kc$ corresponding
to the ideal sheaf \ \m{(\pi_1^{p+1})+(\pi_2^{p+1})}. We have
\begin{eqnarray*}\ko_Z(U_{ij}) & = &
\{(\alpha_1,\alpha_2)\in\ko_{\kc_1}(U_{ij})/(t^{p+1})
\times\ko_{\kc_2}(U_{ij})/(t^{p+1}) \ ; \ \Phi(\alpha_1 \ \text{mod} \ t^p)=
\alpha_2 \ \text{mod} \ t^p \}\\
& = & \{(\alpha_1,\alpha_2)\in\ko_C(U_{ij})[t]/(t^{p+1})\times
\ko_C(U_{ij})[t]/(t^{p+1}) \ ; \ \alpha_1\equiv\alpha_2 \ \text{mod} \ t^p\} .
\end{eqnarray*}
To obtain \m{\ko_{C_2}(U_{ij})}, we have just to quotient by \m{\pi=(t,t)},
and we obtain
\[\ko_{C_2}(U_{ij}) \ = \ \ko_Z(U_{ij})/(t,t) \ \simeq \
 \ko_C(U_{ij})[z]/(z^2) \ , \]
the last isomorphism being
\[(a_0+a_1t+\cdots+a_{p-1}t^{p-1}+\alpha
t^p,a_0+a_1t+\cdots+a_{p-1}t^{p-1}+\beta t^p) \ \mapsto \
\alpha_0+(\beta-\alpha)z . \]
Now we can make explicit the automorphism of \m{\ko_C(U_{ij})[z]/(z^2)} 
induced by
\m{\theta_{ij}} (these isomorphisms will define the cocycle corresponding to
\m{C_2}). It is easy to see that this isomorphism is \m{\phi_{\tau_{ij},1}},
which proves theorem \ref{theo_cf}.
\end{proof}

\newpage

{\bf Note : } This text is a version of 

Dr\'ezet, J.-M. {\em Fragmented deformations of primitive
multiple curves.} Central European Journal of Mathematics 11, n${}^o$ 12 
(2013), 2106-2137.

with the following modifications : 
\begin{enumerate}
\item[--] an unneccessary hypothesis was removed from definition 3.2.2.
\item[--] a supplementary hypothesis was added before proposition 4.5.2.
\end{enumerate}

\end{document}